%% file: main_group_cohomology_hypertoric.tex
\documentclass[a4paper,final]{article}
\usepackage[common_math_textnormal,todonotes,math_base_note,math_simple]{gatmeo}
\usepackage{tikz-cd}

\title{\bf{
The Group Cohomology of Peroidized Hypertoric Variety
}}

\newcommand{\ansonlaw}{\FirstBigRestSmallCaps{Sum kiu Law}\thanks{sklaw@math.cuhk.edu.hk}}
\renewcommand{\omegat}{\FirstBigRestSmallCaps{Nok to Omega Tong}\thanks{onttong@math.cuhk.edu.hk}}
\author{\ansonlaw\ {\small{AND}} \omegat\\\\
\emph{Department of Mathematics, The Chinese University of Hong Kong, Hong Kong}
}
\date{}

\newcommand{\noskipline}{\vspace{-1.5em}}

\renewcommand{\epsilon}{\varepsilon}
\renewcommand{\phi}{\varphi}
\newcommand{\NN}{\mathbb{N}}
\newcommand{\ZZ}{\mathbb{Z}}
\tikzcdset{
    cells={font=\everymath\expandafter{\the\everymath\displaystyle}},
}

\newcommand{\MM}{\mathcal{M}}
\newcommandx*\Mper[1][1=\Gamma]{\mathcal{M}^{\mathrm{per}}(#1)}
\newcommandx*\Mperchar[2][2=\Gamma]{\mathcal{M}^{\mathrm{per}}_{#1}(#2)}
\newcommandx*\Hper[1][1=\Gamma]{\mathcal{A}^{\mathrm{per}}(#1)}
\newcommandx*\Hperchar[2][2=\Gamma]{\mathcal{A}^{\mathrm{per}}_{#1}(#2)}
\newcommandx*\Rper[2][1=\bullet, 2=\Gamma]{\mathcal{R}_{\mathrm{per}}^{#1}(#2)}
\newcommandx*\SRper[2][1=\bullet, 2=\Gamma]{\mathcal{SR}_{\mathrm{per}}^{#1}(#2)}
\newcommandx*\Mfin[1][1=\Gamma]{\mathcal{M}(#1)}
\newcommandx*\Mfinchar[2][2=\Gamma]{\mathcal{M}_{#1}(#2)}
\newcommandx*\Hfinchar[2][2=\Gamma]{\mathcal{A}_{#1}(#2)}
\newcommandx*\Hfin[1][1=\Gamma]{\mathcal{A}(#1)}
\newcommandx*\Rfin[2][1=\bullet, 2=\Gamma]{\mathcal{R}^{#1}(#2)}
\newcommandx*\SR[2][1=\bullet, 2=\Gamma]{\mathcal{SR}^{#1}(#2)}
\newcommandx*\Mmul[1][1=\Gamma]{\mathcal{M}^{\mathrm{mul}}(#1)}
\newcommandx*\idealIper[1][1=\Gamma]{I_{\mathrm{per}}(#1)}
\newcommandx*\idealI[1][1=\Gamma]{I(#1)}
\newcommand{\opendisc}{\mathbb{D}}
\newcommandx*\Betti[1][1=\Gamma]{\mathfrak{B}(#1)}
\newcommandx*\Dolbeault[1][1=\Gamma]{\mathfrak{D}(#1)}

\newcommandx*\corefin[1][1=\Gamma]{\mathcal{L}(#1)}
\newcommandx*\coreper[1][1=\Gamma]{\mathcal{L}^{\mathrm{per}}(#1)}
\newcommandx*\toricP[1][1=P]{\mathcal{X}(#1)}
\newcommandx*\fixptP[1][1=B]{P_{#1}}
\newcommandx*\toricB[1][1=B]{\mathcal{X}(#1)}

\newcommandx*\GC[3][1=\bullet,2=\bullet,3=\Gamma]{\mathfrak{C}^{#1,#2}(#3)} 
\newcommand{\dGC}{\delta} 
\newcommand{\fGC}{\phi}
\newcommand{\gGC}{\psi}
\newcommand{\hGC}{\eta}
\newcommandx*\GCe[4][1=\bullet,2=\bullet,3=\Gamma,4=e]{\mathfrak{C}^{#1,#2}(#3,#4)}

\newcommandx*\GCS[2][1=\bullet,2={\Gamma, S}]{\mathfrak{D}^{#1}(#2)}
\newcommandx*\dGCS[1][1=S]{\dGC_{#1}}
\newcommandx*\fGCS[1][1=S]{\fGC_{#1}}
\newcommandx*\gGCS[1][1=S]{\gGC_{#1}}
\newcommandx*\hGCS[1][1=S]{\hGC_{#1}}
\newcommandx*\CKS[4][1=\bullet,2=\bullet,3=\bullet,4=\Gamma]{\mathfrak{C}^{#1,#2,#3}(#4)} 
\newcommandx*\grCKS[4][1=\bullet,2=\bullet,3=\bullet,4=\Gamma]{I^{#1, #2}\mathfrak{C}^{#3}(#4)} 
\newcommandx*\CKSdiffgrading[3][1=\bullet,2=\bullet,3=\Gamma]{\mathfrak{C}'^{#1,#2}(#3)} 
\newcommand{\dCKS}{d} 
\newcommand{\fCKS}{f}
\newcommand{\gCKS}{g}
\newcommand{\hCKS}{h}
\newcommandx*\HT[3][1=\bullet,2=\bullet,3=\Gamma]{\mathfrak{K}_{\mathrm{per}}^{#1, #2}(#3)}
\newcommandx*\grHT[3][1=\bullet,2=\bullet,3=\Gamma]{I^{#1}\mathfrak{K}_{\mathrm{per}}^{#2}(#3)}
\newcommand{\dHT}{\dCKS_{\mathrm{per}}}
\newcommandx*\HTalt[2][1=\bullet,2=\Gamma]{\tilde{\mathfrak{K}}^{#1}(#2)}

\newcommand{\fHT}{\fCKS_{\mathrm{per}}}
\newcommand{\gHT}{\gCKS_{\mathrm{per}}}
\newcommand{\hHT}{\hCKS_{\mathrm{per}}}
\newcommandx*\HTfin[3][1=\bullet,2=\bullet,3=\Gamma]{\mathfrak{K}^{#1, #2}(#3)}
\newcommandx*\grHTfin[3][1=\bullet,2=\bullet,3=\Gamma]{I^{#1}\mathfrak{K}^{#2}(#3)}
\newcommand{\dHTfin}{\dCKS}
\newcommand{\fHTfin}{\fCKS}
\newcommand{\gHTfin}{\gCKS}
\newcommand{\hHTfin}{\hCKS}
\newcommandx*\CPX[3][1=\bullet,2=\bullet,3=\Gamma]{\mathfrak{D}^{#1,#2}(#3)} 
\newcommandx*\euler[3][3=\Gamma]{\eta^{#1, #2}(#3)}

\newcommandx*\Tutte[2][1=\Gamma]{T(#1;\, #2)}
\newcommandx*\hpoly[2][1=\Gamma]{h(#1;\, #2)}
\newcommandx*\HTutte[2][1=\Gamma]{H(#1;\, #2)}
\newcommandx*\seqn[1][1=n]{[-#1,#1]}
\newcommand{\loopgraph}{{\bigcirc\!\!\bullet}}
\newcommand{\bridgegraph}{{\bullet\mspace{-7mu}-\mspace{-7mu}\bullet}}

\newcommand{\del}{\setminus}
\newcommand{\con}{\mathbin{/}}

\newcommandx*\tree[2][2=\Gamma]{\mathcal{T}_{#2}(#1)}
\newcommandx*\treefunction[1][1=\Gamma]{\mathcal{T}_{#1}}
\newcommandx*\cotreefunction[1][1=\Gamma]{\mathcal{T}_{#1}^*}
\newcommandx*\cotree[2][2=\Gamma]{\mathcal{T}_{#2}^*(#1)}
\newcommandx*\face[2][1=\bullet,2=\Gamma]{F_{#1}(#2)}
\newcommandx*\faceshelling[3][2=\bullet,3=\Gamma]{F_{#2}^{(#1)}(#3)}
\newcommandx*\faceper[2][1=\bullet,2=\Gamma]{F^{\mathrm{per}}_{#1}(#2)}
\newcommandx*\pair[3][3=\Gamma]{\langle #1, #2 \rangle_{#3}}
\newcommandx*\cotreein[2][2=\Gamma]{\mathcal{D}_{#2}(#1)}
\newcommandx*\cotreeinper[2][2=\Gamma]{\mathcal{D}^\per_{#2}(#1)}
\newcommandx*\cotreeinfunction[1][1=\Gamma]{\mathcal{D}_{#1}}
\newcommandx*\externalactivity[2][2=\Gamma]{EA_{#2}(#1)}
\newcommandx*\internalactivity[2][2=\Gamma]{IA_{#2}(#1)}
\newcommand{\intp}[1]{\iota_{#1}} 
\newcommandx*\fundcycle[2][2=\Gamma]{C_{#2}(#1)}
\newcommandx*\fundcyclewithtree[3][2=\Gamma,3=T]{C_{#2}(#3,#1)}
\newcommandx*\basis[2][1=\bullet,2=\Gamma]{\mathcal{B}_{#1}(#2)}
\newcommandx*\basisper[2][1=\bullet,2=\Gamma]{\mathcal{B}^\per_{#1}(#2)}

\newcommandx*\SE[3][1=,3=]{
  \ifstrempty{#3}
  {\ifstrempty{#1}{E^\bullet_{#2}}{#1_{#2}}}
  {\ifstrempty{#1}{E_{#2}(#3)}{#1_{#2}}}
}
\newcommandx*\HMGZ[2][1=\Gamma,2=\bullet]{H^{#2}(\MM(#1),\mathbb{Z})}

\newcommand{\per}{\mathrm{per}}

\newcommand{\HH}{\mathrm{H}}
\newcommandx*\Ext[1][1=\bullet]{\bigwedge\nolimits^{#1}}
\newcommandx*\Sym[1][1=\bullet]{\mathrm{Sym}^{#1}}

\newcommand{\coker}{\operatorname{coker}}
\renewcommand{\im}{\operatorname{Im}}
\newcommand{\supp}{\mathrm{supp}}
\newcommand{\Hom}{\mathrm{Hom}}

\newcommand{\invar}[1]{\Gamma_{#1}}
\newcommandx*\Rder[1][1=\bullet]{R^{#1}}
\newcommandx*\totRder[1][1=\bullet]{\mathbb{R}^{#1}}

\newcommand{\choice}[1]{[#1]} 
\newcommandx*\intord[2][2=S]{[#2]^{#1}} 
\newcommandx*\hyppl[2][2=S]{\alpha^{#2}_{#1}} 

\makeatletter
\newcommand{\GIT}[1][\@nil]{%
  \def\tmp{#1}%
  \ifx\tmp\@nnil
    /\!\!/%
  \else
    /\!\!/_{\! #1}%
  \fi
}
\newcommand{\HQ}[1][\@nil]{%
  \def\tmp{#1}%
  \ifx\tmp\@nnil
    /\!\!/\!\!/\!\!/%
  \else
    /\!\!/\!\!/\!\!/_{\! #1}%
  \fi
}
\makeatother

\begin{document}
\maketitle
\vspace{-3.5em}
\input{abstract}

\thispagestyle{empty}
\tableofcontents
\listoftodos

\input{sec_introduction}

\input{sec_combinatorics}

\input{sec_ht}

\input{sec_periodization_new}
\input{sec_multiplicative_hypertoric_and_cks}
\input{sec_gcnew}
\bibliographystyle{alpha}
\bibliography{ref_group_cohomology_of_periodized_hypertoric}

\end{document}

%% file: abstract.tex
\begin{abstract}
  To a graph $\Gamma$, one can associate a hypertoric variety $\Mfin$ and its multiplicative version $\Mmul$.
  It was shown in \cite{dancso2024deletion} that the cohomology of $\Mmul$ is computed by the CKS complex, which is a finite dimensional complex attached to $\Gamma$.
  The multiplicative hypertoric variety can be realized as the quotient of a periodized hypertoric variety by a lattice action. In this paper, we show that the group cohomology of the lattice with coefficients in the cohomology of the prequotient is isomorphic to the cohomology of the CKS complex using a spectral sequence argument. Therefore, the group cohomology can serve as an alternative way to compute the cohomology of multiplicative hypertoric varieties.
  
  We also found graph-theoretic descriptions for the Euler characteristics of the graded pieces in a certain decomposition of $\HH^\bullet(\Mmul)$.
\end{abstract}

%% file: sec_introduction.tex
\section{Introduction}

Hypertoric varieties are the hyperk\"ahler analogues of toric varieties, defined as the hyperk\"ahler quotient of a complex symplectic vector space by a hyperhamiltonian action of a torus.
The data of this quotient is given by a short exact sequence of tori and a generic character, which is often encoded with a hyperplane arrangment.
In this paper, we focus on hypertoric varieties that come from graphs, and denote by $\Mfin$ the hypertoric variety associated with the graph $\Gamma=(E,V)$. Assuming $\Gamma$ is connected, the associated hyperplane arrangement has $|E|$ hyperplanes in the affine vector space $\mathbb{R}^{|E|-|V|+1}$.

The multiplicative version of hypertoric varieties appears naturally as a special case of multiplicative quiver varieties \cite{Yamakawa_2010} with abelian gauge groups. By analogy with non-abelian Hodge theory, they admit Betti, de Rham and Dolbeault variants, first described by Hausel and Proudfoot \cite{hausel2006untitled}, which have isomorphic cohomology rings. These variants and their relations were further developed in \cite{dancso2024deletion}. The authors showed that the cohomology of the Betti space is isomorphic to the cohomology of the CKS complex $\CKS$ introduced in \cite{cattani19872} and \cite{migliorini2021support}. 

In this paper, we focus on the Dolbeault space, which we denote by $\Mmul$. To define it, we start with the periodized hypertoric variety $\Mper$. This is the hypertoric variety associated with the periodized hyperplane arrangement, namely replace each hyperplane associated to an edge in $\Gamma$ by $\ZZ$ copies of it with integer seperation. Below gives an example of periodization:

\begin{figure}[H]
    \centering
\begin{subfigure}{0.32\textwidth}
\centering
\includegraphics[]{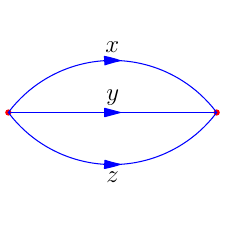}
    \caption{Graph}
    \label{fig:graph_TP2}
\end{subfigure}
\begin{subfigure}{0.32\textwidth}
\centering
\includegraphics[]{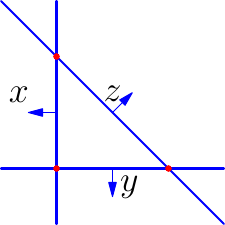}
    \caption{Hyperplane arrangement}
    \label{fig:hyperplane_TP2}
  \end{subfigure}
\begin{subfigure}{0.32\textwidth}
\centering
\includegraphics[]{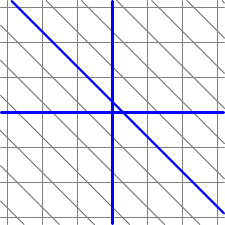}
  \caption{Periodized hyperplane arrangement}
\end{subfigure}
\caption{Example of Graph and Periodization}
\end{figure}

The periodized hyerptoric variety $\Mper$ admits an action by the lattice $\HH_1(\Gamma, \ZZ) \cong \ZZ^{|E|-|V|+1}$. Moreover, the action is free and properly discontinuous on an open retract of $\Mper$. We define $\Mmul$ to be the quotient of the lattice action on this open subspace.

We write $\Rper = \HH^\bullet(\Mper, \ZZ)$ for the cohomology ring of the periodized hypertoric variety. This is related to the cohomology of the corresponding multiplicative hypertoric variety by a Cartan-Leray spectral sequence
\begin{equation*}
    E_2 = \HH^\bullet_{\mathrm{Grp}}(\HH_1(\Gamma, \ZZ), \Rper) \quad \Rightarrow \quad \HH^\bullet(\Mmul, \ZZ),
\end{equation*}
where the $E_2$ page is the group cohomology of the lattice $\HH_1(\Gamma, \ZZ)$ with coefficients in $\Rper$.

We constructed a CKS-type complex $\HT$ together with a homotopy equivalence $ \HT \to \Rper $, which also carries a compatible $\HH_1(\Gamma, \ZZ)$-action. This gives a spectral sequence of double complex
\begin{equation*}
    E_2 = \HH^\bullet(\HH^\bullet_{\mathrm{Grp}}(\HH_1(\Gamma, \ZZ), \HT)) \quad \Rightarrow \quad \HH^\bullet_{\mathrm{Grp}}(\HH_1(\Gamma, \ZZ), \Rper).
\end{equation*}
By some homological algebra arguments, the $E_2$ page is in fact isomorphic to the CKS complex and the spectral sequence degenerates at $E_2$ page, giving us the following main theorem.

\begin{theorem}[thm:]{}
    $\HH_{\mathrm{Grp}}^\bullet(\HH_1(\Gamma), \Rper) \cong \HH^\bullet(\Mmul, \ZZ)$.
\end{theorem}

If we fix an ordering on the edges, the Tutte polynomial of $\Gamma$ is defined as
\[
    \Tutte{x,y}=\sum_{T\text{ spanning tree}}^{}x^{|\internalactivity{T}|}y^{|\externalactivity{T}|},
\]
where $\internalactivity{T}$ and $\externalactivity{T}$ denotes internal and external activity of $T$ respectively. There is an equivalent definition of the Tutte polynomial that does not depend on orders (cref. \Cref{sec:poincare_polynomial_hypertoric}).
The $h$-polynomial of the dual graph of $\Gamma$ is defined to be
\[
    \hpoly{q} := \Tutte{1,q} = 
    \sum_{T\text{ spanning tree}}^{}q^{|\externalactivity{T}|}.
\]

It is proven in \cite{hausel2002toric} and \cite{crawley2004absolutely} that the Poincaré polynomial of $\Mfin$ equals $\hpoly{q}$.
\cite{abdelgadir2022tutte} also gave a proof of this fact by providing a stratification of $\Mfin$ by affine spaces $\mathbb{C}^{d+|\externalactivity{T}|}$ labelled by spanning trees $T$ of $\Gamma$.

The tools we used to construct a homotopy equivalence between $\HTfin$ and $\Rfin$ allow us to construct a monomial basis for $\Rfin$ labelled by spanning trees, where the basis element labelled by $T$ has degree $|\externalactivity{T}|$.
This basis agrees with the one given in \cite[Theorem III.2.5]{stanley1996combinatorics}.
%

We suspect this basis gives a cohomological interpretation of the stratification given in \cite{abdelgadir2022tutte}.

Back to the multiplicative hypertoric varieties, \cite{dancso2024deletion} proved that there is a deletion-contraction sequence of the CKS complex, giving a formula $\chi(\Mmul)=h(1)$ for the Euler characteristic of $\Mmul$, which also equals the number of spanning tree of $\Gamma$.

In this paper, we consider subcomplexes $ \grCKS[k][\ell] = \bigoplus_{p \in \mathbb{N}} \CKS[2p][k-p][\ell] $, giving a decomposition of the CKS complex into a direct sum of these subcomplexes. These pieces inherit the aforementioned deletion-contraction sequence. We showed that their Euler characteristics can be described using a specialization of the Tutte polynomial.

\begin{theorem}{}
    Let $\euler{k}{\ell}$ be the Euler characteristic of the subcomplex $\grCKS[k][\ell]$ of the CKS complex. Then
    \begin{equation*}
        \Tutte{-(x+y+xy),1} = (-1)^d \sum_{k, \ell \in \mathbb{N}} \euler{k}{\ell} \, x^{d-k}y^\ell,
    \end{equation*}
    where $d$ is the genus of $\Gamma$.
\end{theorem}

The contents of this article are organized as follows.
We begin in \Cref{sec:combinatorics} by establishing graph notations. The coherent cotree can be treated as a generalization of external activity that extends coherently to all deletions.

In \Cref{sec:hypertoric_varieties_and_cohomology}, 
we then give a brief introduction of hypertoric varieties and show how $\HTfin$ computes the cohomology ring of $\Mfin$. We also discuss how deletion-contraction of graphs gives rise to different embedding of hypertoric varieties.

In \Cref{sec:periodizaiton}, we extend all the theorems in \Cref{sec:hypertoric_varieties_and_cohomology} to the periodization of graphs and hypertoric varieties.

In \Cref{sec:multiplicative_hypertoric_and_cks}, we explain the $\HH_1(\Gamma,\ZZ)$ action on the periodized hypertoric $\Mper$, the definition of multiplicative hypertoric, and the construction of the CKS complex.

In \Cref{sec:gcnew}, we put everything together and prove the main theorem of this paper. And we finish off with a graph invariant of the multiplicative hypertoric induced by a filtration on the CKS complex.

\nbf{Acknowledgements.}
The authors would like to express the most sincere gratitude to our MPhil advisor Prof. Michael McBreen for his help and guidance. We also thank Ki Fung Chan and Chin Hang Eddie Lam for helpful discussions.

This project was funded by an ECS Grant (reference number 24307121) and GRF Grant (reference number: 14307022) from the Research Grant Council, with PI Michael McBreen.

%% file: sec_combinatorics.tex
\section{Combinatorics}
\label{sec:combinatorics}

\subsection{Graph Notation}

In this paper, a graph means an ``oriented multigraph with loops'' consisting of a set of vertices $V$, a set of edges $E$, and functions $ h, t: E \to V $. We always assume our graphs to be finite and connected, unless stated otherwise. 
In fact, the results of this paper can be generalized to (oriented) matroids. 
However, to keep this paper simple, we will start with graphs and state the properties using the language of matroids.
Readers may refer to \cite{oxley2006matroid} and \cite{bjorner1999oriented} for matroids and oriented matroids.
We will first talk about some properties of graphs ignoring the orientation, i.e. the underlying matroids.

A \vocab{cycle}, or a \vocab{circuit}, of a graph, is a set of edges $Y\subseteq E$, ignoring orientation, that corresponds to a path with the same starting and ending vertices.  
An edge cut of a graph is a subset of edges $X\subseteq E$ such that the number of connected components is increased upon deleting $X$ from the graph. A minimal edge cut will be called a \vocab{bond}, or a \vocab{cocircuit}.

A \vocab{tree} of a graph is a set of edges, ignoring orientation, in which any two vertices are connected by exactly one path.
A \vocab{spanning tree}, or a \vocab{basis}, is a maximal tree. For a connected graph, a spanning tree $T$ has cardinality $|V|-1$ and includes all vertices.
A \vocab{spanning cotree}, or a \vocab{cobasis}, is the complement of a spanning tree in the set of edges. A \vocab{cotree} is a subset of a spanning cotree.

Specifying the set of all cycles allows one to describe all the spanning trees: the maximal subsets of $E$ that do not contain any cycle. The converse is also true. The same relation holds for bonds and cotrees. Indeed, the correspondence between bases and cobases arises from a duality of matroids. We have the following propositions.

\begin{proposition}[pps:]{}
    The collection of all cycles and the collection of all bonds both satisfy the (weak) circuit elimination axiom:

    Let $C_1$ and $C_2$ be distinct elements in the collection. If $e\in C_1\cap C_2$, then the collection contains a member $C_3$ such that $C_3\subseteq C_1\cup C_2\del e$.
\end{proposition}

\begin{proposition}[pps:]{}
    The collection of all spanning trees and the collection of all spanning cotrees both satisfy the symmetric basis exchange axiom:

    Let $A$ and $B$ be distinct elements in the collection. If $a\in A\del B$, then there exists an element $b\in B\del A$ such that both $(A\del a)\cup b$ and $(B\del b)\cup a$ are in the collection.
\end{proposition}

This property guarantees that all spanning trees and spanning cotrees have the same cardinality. In a connected graph, the cardinality of a spanning tree and spanning cotree is $|V|-1$ and $|E|-|V|+1$ respectively. We call $d=g(\Gamma):=|E|-|V|+1$ the \vocab{genus} of the graph, the usual definition of genus if you view the graph as a simplicial complex.

As described above, an equivalent definition for cycle is the minimal set of edges that are not contained in a spanning tree, using the fact that all spanning trees have cardinality $|V|-1$, and any $|V|$ edges contain a cycle.
Dually, any $d+1$ distinct edges in a genus $d$ graph contain a bond.

To a graph, we can associate with it a simplicial complex on the set $E$, whose faces are subsets of $E$ contained in some spanning tree, called the independence complex.
In this paper, we focus more on the dual structure of the graph.

\begin{definition}[def:face_coindependent_complex]{}
    We define $\face$ as the coindependence complex whose faces are subsets of $E$ that are contained in some spanning cotree, or equivalently, the subsets of $E$ that does not contain a bond.
    We let $\face[k]$ be the set of cardinality $k$ elements of $\face$. In particular, 
    \[F_k(\Gamma) := \{ S \subseteq E \mid |S| = k \text{ and } S \text{ does not contain a bond} \}.\]
\end{definition}

\subsection{Graph Homology and Cohomology}
\label{sec:graph_homology_and_cohomology}

Any graph can be enhanced to a 1-truncated simplicial set (equivalently, reflexive graph) by adding a degenerated edge at each vertex. By a morphism of graphs, we mean a morphism of the associated simplicial set. Therefore, graphs inherit the notion of homology and cohomology of the associated simplicial sets, which are described as follows.

Let $A$ be an abelian group and define $ \partial_{\Gamma, A} : A^E \to A^V $ by $ e \mapsto h(e) - t(e) $ for all $ e \in E $, and let $ d_{\Gamma, A} : A^V \to A^E $ be its transpose.

\begin{definition}[def:]{}
    The graph homology and cohomology of $\Gamma$ with coefficients in $A$ are defined as follows.
    \begin{align*}
        \HH_0(\Gamma, A) &= \coker(\partial_{\Gamma, A}) & \HH_1(\Gamma, A) &= \ker(\partial_{\Gamma, A}) \\
        \HH^0(\Gamma, A) &= \ker(d_{\Gamma, A}) & \HH^1(\Gamma, A) &= \coker(d_{\Gamma, A})
    \end{align*}
    \noskipline
    \begin{remark}
        If unspecified, graph homology and cohomology take coefficients in $\ZZ$ by default.
    \end{remark}
    \begin{remark}
        We always have $ \HH_0(\Gamma, A) \cong \HH^0(\Gamma, A) \cong A^{\pi_0(\Gamma)} $ canonically, where $\pi_0(\Gamma)$ is the set of connected components of $\Gamma$. Non-canonically, we also have $ \HH_1(\Gamma, A) \cong \HH^1(\Gamma, A) \cong A^{g(\Gamma)} $, where $g(\Gamma)$ is the genus of $\Gamma$.
    \end{remark}
    \begin{remark}
        Notice that if $A$ is a field, $\partial_\Gamma$ is a representation of $\Gamma$, i.e. if we identify $e\in E$ with its basis vector in $A^E$, then $S\subseteq E$ is contained in a spanning tree if and only if $\partial_\Gamma|_{S}$ is injective and $\partial_\Gamma(S)$ is linearly independent. And by matroid theory, the natural projection $\pi_{\Gamma}:A^E\rightarrow \HH^1(\Gamma)$ is a dual representation (replace spanning tree with spanning cotree in the above definition) of $\Gamma$.
    \end{remark}
\end{definition}

The linear combination of edges (with orientation) in a cycle is in the kernel of $\partial_\Gamma$, hence giving an element in $\HH_1(\Gamma)$. On the other hand, any edge $e \in E$ gives a class $ [e] \in \HH^1(\Gamma) $ through the projection $ \ZZ^E \to \HH^1(\Gamma) $. 

\begin{definition}[def:pair_homology_edge]{}
    The transpose of the inclusion $ \HH_1(\Gamma) \hookrightarrow \ZZ^E $ gives a pairing $ \pair{-}{-} : \HH_1(\Gamma) \otimes \ZZ^E \to \ZZ $. Given $ S \subseteq E $, we define $ \pair{-}{S} : \HH_1(\Gamma) \to \ZZ^S $ by $ \pair{\gamma}{S} = (\pair{\gamma}{e})_{e \in S} $.
\end{definition}

In this paper, we often need to fix a basis for $\HH_1(\Gamma)$ and $\HH^1(\Gamma)$. One way to do it is to pick a spanning tree (or equivalently a spanning cotree).

\begin{definition}[def:fundamental_cycle]{}
    Upon choosing a spanning tree $T$, any edge $e\in E\setminus T=:T^*$, the corresponding spanning cotree, determines a unique cycle $\fundcyclewithtree{e}\subseteq T\cup e$ that contains $e$, called the \vocab{fundamental cycle}.
    We will abbreviate $\fundcyclewithtree{e}$ as $\fundcycle{e}$ if the choice of spanning tree is fixed.
\end{definition}

By abuse of notation, we denote by $\fundcyclewithtree{e}$ also the corresponding element in $\HH_1(\Gamma)$ whose coefficient of $e$ is $1$. This extends to a linear map $ \fundcyclewithtree{-} : \ZZ^{T^*} \to \HH_1(\Gamma) $.


\begin{theorem}[thm:fundamental_cycle_basis_for_graph_homology]{}
    Fixing a spanning tree $T \subseteq E$, the map $ \fundcycle{-}: \ZZ^{T^*} \to \HH_1(\Gamma) $ is an isomorphism, whose inverse is given by the composite $ \HH_1(\Gamma) \hookrightarrow \ZZ^E \twoheadrightarrow \ZZ^{T^*} $.
    \begin{remark}
        In particular, $\{\fundcycle{x}\}_{x \in T^*}$ forms a basis of $\HH_1(\Gamma)$. Dually, $\{[x]\}_{x \in T^*}$ forms a basis of $\HH^1(\Gamma)$.
    \end{remark}
    \begin{proof}
        For all $ x, y \in T^* $, $\pair{\fundcycle{x}}{y}$ is by definition $1$ if $x=y$, and $0$ else. Thus, the composite $ \ZZ^{T^*} \to \HH_1(\Gamma) \hookrightarrow \ZZ^E \twoheadrightarrow \ZZ^{T^*} $ is the identity. On the other hand, for all $ \gamma \in \HH_1(\Gamma) $, $\gamma-\sum_{x\in T^*}\pair{\gamma}{x}\fundcycle{x}$ must be $0$ since it is in $\ZZ^T$, which contains no cycles. Thus, the other composite is also the identity.
    \end{proof}
\end{theorem}

Under the identification $ \HH_1(\Gamma) \cong \ZZ^{T^*} $, the pairing
in \Cref{def:pair_homology_edge} gives a pairing $ \ZZ^{T^*} \otimes \ZZ^E \to \ZZ $, which we also denote by $\pair{-}{-}$. One can check that $\pair{x}{e}\in \CB{-1,0,1}$ for all $x\in T^*$ and $e\in E$. Moreover, the matrix defining the pairing is totally unimodular, i.e. every square submatrix has a determinant in $\CB{-1,0,1}$.

\subsection{Deletion and Contraction of Graphs}
\label{sec:deletion_contraction_graph_theory}

Let $ e \in E $. The deletion $\Gamma \del e$ is defined to be the subgraph of $\Gamma$ removing the edge $e$, which has vertices $V$ and edges $E \setminus e$. The contraction $\Gamma \con e$ is defined to be the quotient of $\Gamma$ contracting the edge $e$, which has vertices $V/(h(e)=t(e))$ and edges $E \setminus e$. One can show that the deletion and contraction of multiple edges do not depend on order. Hence we write $\Gamma\del\CB{x_1,\dots ,x_n}=\Gamma\del x_1\del \cdots \del x_n$ and $\Gamma\con\CB{x_1,\dots ,x_n}=\Gamma\con x_1\con\cdots \con x_n$ for $x_1,\dots ,x_n\in E$ distinct.

For any $ S \subseteq E $, deletion and contraction give maps $ \Gamma \del S \to \Gamma $ and $ \Gamma \to \Gamma \con S $ respectively. If $S$ is a tree, the maps induced by $ \Gamma \to \Gamma \con S $ are isomorphisms. If $S$ is a cotree, we have a short exact sequence
\begin{equation} \label{eqn:del_seq}
  \begin{tikzcd}
      0 \arrow[r] & \HH_1(\Gamma \del S) \arrow[r] & \HH_1(\Gamma) \arrow[r] & \ZZ^S \arrow[r] & 0
  \end{tikzcd}
\end{equation}
where the projection is the composite $ \HH_1(\Gamma) \hookrightarrow \ZZ^E \twoheadrightarrow \ZZ^S $. In the case where $ S = \{ e \} $ is not a bond, the following short exact sequence defines the ``interior product'' $\intp{e}$ with the class $ [e] \in \HH^1(\Gamma) $.
\begin{equation} \label{eqn:intp}
    \begin{tikzcd}
        0 \arrow[r] & \Ext \HH_1(\Gamma \del e) \arrow[r] & \Ext \HH_1(\Gamma) \arrow[r, "\intp{e}"] & \Ext[\bullet-1] \HH_1(\Gamma \del e) \arrow[r] & 0
    \end{tikzcd}
\end{equation}

If $e$ is not a bond, $ S \mapsto S \cup e $ defines a map $ \face[\bullet-1][\Gamma \del e] \to \face $, whose image is the set of $ S \in \face $ such that $ e \in S $. On the other hand, $ \face[\bullet][\Gamma \con e] = \{ S \in \face \mid e \notin S \} $. This gives the following short exact sequence.
\begin{equation} \label{eqn:face_delcon}
    \begin{tikzcd}
        0 \arrow[r] & \ZZ^{\face[\bullet-1][\Gamma \del e]} \arrow[r] & \ZZ^{\face} \arrow[r] & \ZZ^{\face[\bullet][\Gamma \con e]} \arrow[r] & 0
    \end{tikzcd}
\end{equation}

\missing[ok]{rewrite}
%

\subsection{Coherent Tree}
\label{sec:coherent_tree}

Considering all deletions at once, there is a direct system of abelian groups $\{\HH_1(\Gamma \del S)\}_{S \in \face}$. Picking a spanning tree for each deletion gives a basis for all these abelian groups. We can make such choices ``coherently'', giving an inverse system of splits to the direct system, so that each split is a coordinate projection. We will see that these choices relate to the external activities in the graph.

\subsubsection{Coherent (Co)tree and Graph (Co)homology}

\begin{definition}[def:coherent_tree]{}
  Let $\treefunction$ (resp. $\cotreefunction$)  be a function from $\face$ to $2^E$. We say $\treefunction$ (resp. $\cotreefunction$) is a coherent (co)tree if it satisfies the following:
  \begin{enumerate}
    \item $\tree{S}$ (resp. $\cotree{S}$) is a spanning (co)tree in $\Gamma\del S$.\label{cond:coherent_tree_1}
    \item If $S_2\subseteq S_1\subseteq E$, then $\tree{S_2}\cup S_2\subseteq \tree{S_1}\cup S_1$ (resp. $\cotree{S_1}\subseteq \cotree{S_2}$).\label{cond:coherent_tree_2}
  \end{enumerate}
  \begin{remark}
    There is a bijection between coherent trees and cotrees: any coherent tree will give a coherent cotree by taking complement inside $\Gamma\del S$ for each $S$, and vice versa.
  \end{remark}
  \begin{remark}
    Since the domain of $\treefunction$ is $\face$, the graph $\Gamma\del S$ is always connected. Hence for condition 1, we are guaranteed to be able to find a spanning tree (instead of a spanning forest). And the cardinality $|\tree{S}|=|V|-1$ is hence constant for $S\in\face$. Dually, the cardinality $|\cotree{S}|=|E|-|V|+1-|S|=d-|S|$.
  \end{remark}
  \begin{remark}
    A coherent (co)tree of $\Gamma$ determines a coherent (co)tree of $\Gamma \del S$ for all $ S \in \face $ by setting $ \tree{S'}[\Gamma \del S] = \tree{S \cup S'} $ (resp. $ \cotree{S'}[\Gamma \del S] = \cotree{S \cup S'} $).
  \end{remark}
  \begin{remark}
    In general, the construction of coherent (co)tree is not compatible with contraction of $\Gamma$. However, if $t\in \tree{\emptyset}$, we can extend the construction of a coherent cotree to $\Gamma \con t$.
    In particular, we can set $\cotree{S}[\Gamma\con t] = \cotree{S}$ for all $S\in \face[\bullet][\Gamma\con t]$.
  \end{remark}
\end{definition}

A coherent (co)tree satisfies the following property.
\begin{lemma}[lem:tree_invariant]{}
    Let $\treefunction$ (resp. $\cotreefunction$) be a coherent (co)tree of $\Gamma$ and let $ S \in \face $. For all $ x \in \cotree{S} $, we have $ \tree{S \cup x} = \tree{S} $ (resp. $ \cotree{S \cup x} = \cotree{S} \setminus x $).
    \begin{proof}
        By definition, we have $ \tree{S} \subseteq \tree{S \cup x} \cup x $. Since $ x \notin \tree{S} $ and $ |\tree{S \cup x}| = |\tree{S}| $, it follows that $ \tree{S \cup x} = \tree{S} $.
    \end{proof}
\end{lemma}

To prove the existence of a coherent (co)tree, we need to introduce the notion of \vocab{shelling}. Let $\Delta$ be a pure simplicial complex. Fix a linear ordering $T_1,T_2,\dots $ of maximal simplices of $\Delta$. Let $\Delta^{(k)}$ be the simplicial complex generated by the first $k$ maximal simplicies, i.e. $\Delta^{(k)}=\CB{G\in \Delta\mid G\subseteq T_i, i\leq k}$. We say such a linear ordering is a \vocab{shelling} if the complex $\Delta^{(k-1)} \cap T_k$ is pure and of dimension $\dim T_k-1$ for all $k=2,3,\dots $. A simplicial complex is said to be \vocab{shellable} if it admits a shelling. Let $R(T_k)=\CB{x\in T_k\mid T_k\del x\in \Delta^{(k-1)}}$ be the \vocab{restriction} of $T_k$ induced by the shelling. We have the following classical result:
\begin{proposition}[pps:shelling_RT_unique_minimal]{\cite[Proposition 7.2.2]{bjorner1992homology}}
  Given a shelling $T_1,T_2,\dots $ of $\Delta$, $R(T_k)$ is the unique minimal subset of $T_k$, with respect to inclusion, that does not lie in $\Delta^{(k-1)}$. In particular,
  \[
  \Delta^{(k)}=\Delta^{(k-1)}\cup \CB{S\mid R(T_k)\subseteq S\subseteq T_k}.
  \]
\end{proposition}

The coindependence complex $\face$ is a pure simplicial complex since its maximal simplices are spanning cotrees, which have the same cardinality. After fixing an ordering on the set $E$, it induces a linear order, the lexicographic order, on the spanning cotrees. It is proven in \cite{bjorner1992homology} that such a linear order on the spanning cotrees is a shelling. In general, the (co)independence complex of a matroid is shellable.

We will inductively construct a coherent cotree using the shelling induced by a fixed ordering on the set $E$. Let $T^*_1, T^*_2,\dots $ be the shelling. For all subsets $S\subseteq T^*_1$, we let $\cotree{S}=T^*_1\del S$. This defines $\cotreefunction$ on $\faceshelling{1}$. To extend $\cotreefunction$ from $\faceshelling{k-1}$ to $\faceshelling{k}$, by \Cref{pps:shelling_RT_unique_minimal}, we only need to define $\cotreefunction$ for $\CB{S\mid R(T^*_k)\subseteq S\subseteq T^*_k}$. We let $\cotree{S}=T^*_k\del S$ for all such $S$. We can check that this construction gives us a coherent cotree:

\begin{theorem}[thm:]{}
  For any graph $\Gamma$, there exists a coherent (co)tree.
  \begin{proof}
    The construction above satisfies Condition \ref{cond:coherent_tree_1} of \Cref{def:coherent_tree}. We will prove that $\cotreefunction$ satisfies Condition \ref{cond:coherent_tree_2} on $\faceshelling{k}$ inductively. It is obviously true for $k=1$. Assume it is true for $k-1$. To prove that it is true for $k$, we need to check Condition \ref{cond:coherent_tree_2}. Let $ S_1, S_2 \subseteq T_k^* $ so that at least one of them does not lie in $\faceshelling{k-1}$. If $S_1, S_2 \notin \faceshelling{k-1}$, by construction $\cotreefunction$ satisfies Condition \ref{cond:coherent_tree_2} as well. The remaining case is that $S_2$ is in $\faceshelling{k-1}$ but $S_1$ is not. In fact, we just need to prove that $\cotree{R(T^*_k)\del e}\supseteq \cotree{R(T^*_k)}$ for all $e\in R(T^*_k)$ as any $S\subsetneq R(T^*_k)$ is a subset of $R(T^*_k)\del e$ for some $e$ and by induction $\cotree{S}\supseteq \cotree{R(T^*_k\del e)}$.

    Let $k'<k$ be the smallest index such that $T^*_{k'}\supseteq R(T^*_k)\del e$. By construction, $\cotree{R(T^*_k)\del e}=T^*_i\del(R(T^*_k)\del e)$. By the definition of $R(T^*_k)$, $T^*_k\del e\in \faceshelling{k-1}$, hence there exist some $i<k$ such that $T^*_i\supseteq T^*_k\del e$, i.e. $T^*_i=T^*_k\del e\cup f$ for some $f<e$. We claim that $k'$ equals the smallest index $i$ such that $T^*_i\supseteq T^*_k\del e$. Indeed, if this is the case, then $\cotree{R(T^*_{k'})\del e}=T^*_{k'}\del(R(T^*_k)\del e)=T^*_k\cup f\del R(T^*_k)\supseteq T^*_k\del R(T^*_k)=\cotree{R(T^*_k)}$, which concludes the proof.

    Assume contrary that $k'<i$, we have $R(T^*_k)\del e\subseteq T^*_{k'},T^*_i$. Let $h\in T^*_{k'}\del(R(T^*_k)\del e)$ be the smallest element not in $T^*_i$, by the symmetric basis exchange axiom of cotrees, there exists an element $g\in T^*_i\del(R(T^*_k)\del e)$ such that $T^*_j=T^*_i\del g\cup h$ is a cotree. Note that $k'<i$ and $h<g$. We will construct a cotree with a smaller index that also contains $R(T^*_k)$, which raises a contradiction.
    \begin{itemize}
      \item If $T^*_k\del g\cup h\supseteq R(T^*_k)$ is a cotree, then we are done.
      \item Else, $T^*_k\del g\cup h$ contains a bond $B$, which necessary contains $h$ since otherwise $B\subseteq T^*_k\del h\subseteq T^*_k$, contradicting $T^*_k$ is a cotree. A similar argument shows that $B$ contains $e$ using the fact that $T^*_j$ is a cotree. Since $T^*_k$ is a cotree, $T^*_k\cup f$ contains a bond $D$ which contains $f$. By a similar argument, $D$ contains $e$ as $T^*_i=T^*_k\del e\cup f$ is a cotree.
        \begin{itemize}
          \item If $D$ does not contain $g$, notice that $C\neq D$ as $f\notin C$ but $f\in D$. By the circuit elimination axioms on the bonds, there is a bond contained in $C\cup D\del e\subseteq T^*_k\del e\del g\subseteq T^*_j$, contradicting $T^*_j$ is a cotree.
          \item Else if $f<g$, We claim that $T^*_k\del g\cup f$ is a basis. Assume otherwise that it contains a bond $D'$, which as above contains $f$. Since $g\in D$ but $g\notin D'$, $D\neq D'$, and we can use the circuit elimination axioms on $D'$ and $D$ to imply that there is a bond contain in $D'\cup D\del f\subseteq T^*_k$, contradiction.
          \item Else $f\geq g$, in particular we have $e>f>h$. Using the fact that $e,h\in C$, $T^*_k\del e\cup h$ is a basis. However, $T^*_k\del e\cup h$ is smaller than $T^*_i=T^*_k\del e\cup f$, contradicting the assumuption of $T^*_i$.
        \end{itemize}
    \end{itemize}
  \end{proof}
\end{theorem}

\begin{example}[exp:coherent_cotree_TP2]{}
     The graph $\Gamma$ in \Cref{fig:graph_TP2} has a bond $\CB{x,y,z}$. $\face=2^E\del \CB{x,y,z}$. The shelling induced by the ordering $x<y<z$ on $E$ is $T^*_1=\CB{x,y}, T^*_2=\CB{x,z}$ and $T^*_3=\CB{y,z}$. $\faceshelling{1}=\CB{\emptyset ,\CB{x},\CB{y},\CB{x,y}}$. 
    Since $\CB{z}=T^*_2\del x\notin \faceshelling{1}$ and $\CB{x}=T^*_2\del z\in \faceshelling{1}$, $R(T^*_2)=\CB{z}$.
    \Cref{pps:shelling_RT_unique_minimal} says that 
    \[
        \faceshelling{2}=\CB{\emptyset ,\CB{x},\CB{y},\CB{z},\CB{x,y},\CB{x,z}}=\faceshelling{1}\cup \CB{\CB{z},\CB{x,z}}.
    \]
    The coherent cotree induced from this shelling is defined by
    \[
        \begin{array}{c|ccccccc}
            S\in \face & \emptyset & \CB{x} & \CB{y} & \CB{z} & \CB{x,y} & \CB{x,z} & \CB{y,z}\\
            \hline
            \cotree{S} & \CB{x,y} & \CB{y} & \CB{x} & \CB{x} & \emptyset & \emptyset & \emptyset
        \end{array}
    \]
\end{example}


\subsubsection{Coordinate Based Map between Graph (Co)homologies}

Using \Cref{thm:fundamental_cycle_basis_for_graph_homology}, a coherent cotree $\cotreefunction$ determines isomorphisms $ \HH_1(\Gamma \del S) \cong \ZZ^{\cotree{S}} \cong \HH^1(\Gamma \del S) $. And the inclusion $ \HH_1(\Gamma \del S) \hookrightarrow \HH_1(\Gamma) $ can be expressed in coordinates as
\begin{equation} \label{eqn:hom_incl_in_coord}
    x \mapsto x + \sum_{y \in \cotree{\emptyset} \setminus \cotree{S}} \pair{x}{y}[\Gamma \del S] y,
\end{equation}
where we extended the domain of the pairing $ \pair{-}{-}[\Gamma \del S] : \ZZ^{\cotree{S}} \otimes \ZZ^{E \setminus S} \to \ZZ $ to $ \ZZ^{\cotree{S}} \otimes \ZZ^E $ along the projection $ \ZZ^E \twoheadrightarrow \ZZ^{E \setminus S} $. From this formula, it is clear that the coordinate projection $ \ZZ^{\cotree{\emptyset}} \twoheadrightarrow \ZZ^{\cotree{S}} $ splits the inclusion $ \HH_1(\Gamma \del S) \hookrightarrow \HH_1(\Gamma) $.

From \Cref{eqn:del_seq} and \Cref{eqn:hom_incl_in_coord}, it is clear that $\{ \pair{x}{S} \}_{x \in \cotree{\emptyset} \setminus \cotree{S}}$ is a basis of $\ZZ^S$. Denote by $\{ \hyppl{x} \}_{x \in \cotree{\emptyset} \setminus \cotree{S}}$ the dual basis, we have the following:

\begin{lemma}[lem:pairing_for_del]{}
    $ \pair{x}{y}[\Gamma \del S] = -\hyppl{y}(\pair{x}{S}) $ for all $ x \in \cotree{S} $ and $ y \in \cotree{\emptyset} \setminus \cotree{S} $.
    \begin{proof}
        \Cref{eqn:hom_incl_in_coord} implies that
        \begin{equation*}
            \fundcycle{x}[\Gamma \del S] = \fundcycle{x} + \sum_{y \in \cotree{\emptyset} \setminus \cotree{S}} \pair{x}{y}[\Gamma \del S] \fundcycle{y}
        \end{equation*}
        in $\HH_1(\Gamma)$ for all $ x \in \cotree{S} $. Applying $ \pair{-}{S} : \HH_1(\Gamma) \to \ZZ^S $ to the above gives
        \begin{equation*}
            \pair{x}{S} = - \sum_{y \in \cotree{\emptyset} \setminus \cotree{S}} \pair{x}{y}[\Gamma \del S] \pair{y}{S}.
        \end{equation*}
        Thus, $ \pair{x}{y}[\Gamma \del S] = -\hyppl{y}(\pair{x}{S}) $.
    \end{proof}
\end{lemma}
\missing[a]{wanna rewrite this later}

\subsubsection{A Subset Induced by Coherent Cotree and External Activity}

We define the following operator induced by coherent cotree, which may roughly be thought of as an equivalent definition of coherent cotree. But we will not prove it here. The main purpose of this operator is for us to find a basis $\basis$ for $\Rfin$ (cref. \Cref{def:Rfin}) as a vector space. See \Cref{thm:Rfin_split}.

\begin{definition}[def:]{}
  Fix a coherent cotree $\cotreefunction$ in $\Gamma$. For each $ S \in \face $, define
  \begin{equation*}
    \cotreein{S} = \{ e \in S \mid e \in \cotree{S \setminus e} \}=\CB{e\in S\mid \cotree{S\setminus e}=\cotree{S}\cup e}.
  \end{equation*}
  And let $ \basis = \{ S \in \face \mid \cotreein{S} = \emptyset \} $.
\end{definition}

\begin{example}[exp:cotreein_TP2]{}
    Continuing \Cref{exp:coherent_cotree_TP2}, we have
\[
    \begin{array}{c|ccccccc}
        S\in \face & \emptyset & \CB{x} & \CB{y} & \CB{z} & \CB{x,y} & \CB{x,z} & \CB{y,z}\\
        \hline
        \cotreein{S} & \emptyset & \CB{x} & \CB{y} & \emptyset & \CB{x,y} & \CB{x} & \emptyset
    \end{array}
\]
Hence $ \basis = \CB{\emptyset, \CB{z}, \CB{y,z}}$.
\end{example}

A familiar concept in graph theory is external activity.

\begin{definition}[def:]{}
  Fix an ordering on $E$. And let $T$ be a spanning tree in $\Gamma$. We say an edge $e\in E\del T$ is externally active with respect to $T$ if $e$ is the smallest edge in $\fundcyclewithtree{e}$ with respect to the ordering on $E$. And denote $\externalactivity{T}$ to be the set of externally active edges with respect to $T$.
\end{definition}

The two definitions are related by the following proposition.

\begin{proposition}[pps:coherent_cotree_and_external_activity]{}
  Let $T^*$ be a spanning cotree in $\Gamma$. If the coherent cotree is constructed from an ordering on the set $E$ and external activity is defined with the same ordering, then $\cotreein{T^*}=\externalactivity{T=E\del T^*}$.
  \begin{proof}
    Note that $\cotreein{T^*\del e}$ is a singleton. 
    By the construction of $\cotreefunction$, 
    \[
      \cotreein{T^*}=\CB{e\in T^*\mid \CB{e}=\cotree{T^*\setminus e}}=\CB{e\in T^*\mid T^*\textnormal{ is the smallest spanning cotree containing }T^*\del e},
    \]
    where the ordering is with respect to the shelling induced by the ordering on $E$.
    Consider the flat $F$ of $T^*\del e$, i.e. the largest subset of $E$ containing $T^*\del e$ that does not contain a spanning cotree. By matroid theory, $E-F=\fundcyclewithtree{e}$. Hence $e$ is the smallest element in $\fundcyclewithtree{e}$ if and only if $T^*$ is the smallest spanning cotree containing $T^*\del e$.
  \end{proof}
\end{proposition}

For contraction $ \Gamma \con t $ with $ t \in \tree{\emptyset}$, which is given a coherent cotree as described in \Cref{def:coherent_tree}, we have $ \cotreein{S}[\Gamma \con t] = \cotreein{S} $ for all $ S \in \face[\bullet][\Gamma \con t] $.

For deletion $\Gamma\del S'$ with $S'\subseteq \face$, we have
\[
  \cotreein{S}[\Gamma\del S']=\{ e \in S \mid e \in \cotree{S \setminus e}[\Gamma\del S'] \}=\{ e \in S \mid e \in \cotree{S\cup S' \setminus e} \}= \cotreein{S\cup S'}\cap S
\]
for all $ S \in \face[\bullet][\Gamma\del S'] $.
We also have the following lemma about $\cotreein{S}$, which resembles \Cref{lem:tree_invariant}.

\begin{lemma}[lem:cotreein_subset]{}
  Let $ S \in \face $. Let $e\in E\del S$ such that $S\cup e\in \face$. Then $\cotreein{S}[\Gamma\del e]\subseteq \cotreein{S}$. We also have
  \begin{itemize}
    \item $ \cotreein{S \cup t} = \cotreein{S}[\Gamma\del t] \subseteq \cotreein{S} $ if $ t \in \tree{S} $, and
    \item $ \cotreein{S \cup x} = \cotreein{S}[\Gamma\del x]\cup x= \cotreein{S} \cup x $ if $ x \in \cotree{S} $.
  \end{itemize}
  \begin{proof}
    They both follow directly from the definition of coherent cotree. For the first statement,
    \[
      y\in \cotreein{S}[\Gamma\del e]\iff y\in \cotree{S\del y}[\Gamma\del e]=\cotree{S\del y\cup e}\subseteq \cotree{S\del y},
    \]
    which implies $y\in \cotreein{S}$. 

    If $t=e\in \tree{S}$, then $t\notin \cotreein{S\cup t}$, i.e. $\cotreein{S\cup t}=\cotreein{S\cup t}\cap S=\cotreein{S}[\Gamma\del t]$.
    If we have $x=e\in \cotree{S}$, by \Cref{lem:tree_invariant}, $\cotree{S}=\cotree{S\cup x}\cup x$, hence $x\in \cotreein{S\cup x}$. The first equality hence follows similarly. For the second equality, we need to prove the reverse inclusion $\cotreein{S}\subseteq \cotreein{S}[\Gamma\del x]=\cotreein{S\cup x}\cap S$. If $e\in \cotreein{S}$, i.e. $\cotree{S\del e}=\cotree{S}\cup e$, then $\cotree{S\del e}=\cotree{S\cup x}\cup x \cup e$. Consider $S\del e\cup x$, we have 
    \[
      \cotree{S\del e}\subseteq \cotree{S\del e\cup x}\subseteq \cotree{S\cup x}=\cotree{S\del e}\cup e \cup x.
    \]
    Since $x\notin \cotree{S\del e\cup x}$, by the definition of coherent cotree $\cotree{S\del e\cup x}=\cotree{S\cup x}\cup e$, $e\in \cotreein{S\cup x}$. Hence $\cotreein{S}\subseteq \cotreein{S\cup x}\cap S$.
  \end{proof}
\end{lemma}

\begin{corollary}[crl:basis_and_cotreein]{}
  $\basis=\CB{T^*\del \cotreein{T^*}\mid T^* \subseteq E \textnormal{ is a spanning cotree}}$.
  \begin{proof}
    By \Cref{lem:cotreein_subset}, $T^*\del \cotreein{T^*}\in \basis$ for all spanning cotree $T^*$. The inverse map is given by observing $S\cup\cotree{S}$ is a spanning cotree.
  \end{proof}
\end{corollary}

%% file: sec_ht.tex
\section{Hypertoric Varieties and their Cohomology}
\label{sec:hypertoric_varieties_and_cohomology}

\subsection{Hypertoric Varieties}
\label{sec:hypertoric_varieties}

We now give a short introduction to hypertoric varieties. A more extensive account could be found in \cite{proudfoot2007survey}.

Given a short exact sequence of complex tori
\begin{equation*}
    \begin{tikzcd}
        1 \arrow[r] & G \arrow[r, "\iota"] & D \arrow[r, "\pi"] & T \arrow[r] & 1
    \end{tikzcd}
\end{equation*}
with $ D \cong (\mathbb{C}^\times)^n $, we consider
\begin{enumerate}
    \item the $G$-action on $T^*\mathbb{C}^n$ by $ g \cdot (z, w) = (\iota(g)_iz_i, \iota(g)_i^{-1}w_i)_{i=1}^n $, and
    \item the moment map $ \mu : T^*\mathbb{C}^n \to \mathfrak{g}^* $ for the $G$-action given by $ (z, w) \mapsto \iota^*((z_iw_i)_{i=1}^n) $.
\end{enumerate}
Together with a character $ \chi : G \to \mathbb{C}^\times $ and $ \lambda \in \mathfrak{g}^* $, we define the hypertoric variety associated with these data to be the GIT quotient
\begin{equation*}
    \mathcal{M}_{\chi, \lambda} = T^*\mathbb{C}^n \HQ[\chi, \lambda] G = \mu^{-1}(\lambda) \GIT[\chi] G.
\end{equation*}
There is a remaining $T$-action on the quotient with a moment map $ \nu : \mathcal{M}_{\chi, \lambda} \to \mathfrak{t}^* $. We will always set $\lambda=0$. In this case, the above data can be assembled into a hyperplane arrangement with $n$ hyperplanes in the affine subspace $ (\iota_{\mathbb{R}}^*)^{-1}(\chi) = \{ x \in \mathfrak{d}_{\mathbb{R}}^* \mid \iota_{\Gamma}^*(x) = \chi \} $ of $\mathfrak{d}_{\mathbb{R}}^*$, where the $i$-th hyperplane is the intersection of $(\iota_{\mathbb{R}}^*)^{-1}(\chi)$ with the $i$-th coordinate hyperplane in $\mathfrak{d}_{\mathbb{R}}^*$. Conventionally, one choose a lift $\theta$ of $\chi$ along $\iota_\mathbb{R}^*$, allowing us to identify the affine subspace $(\iota_{\mathbb{R}}^*)^{-1}(\chi)$ with $\mathfrak{t}_\mathbb{R}^*$. This gives a hyperplane arrangement
\begin{equation*}
    \mathcal{A}_\chi = \{ H_i = \{ x \in \mathfrak{t}_{\mathbb{R}}^* \mid \pi^*(x)_i = \theta_i \} \}_{i=1}^n.
\end{equation*}
Different choices of lifts give the same hyperplane arrangement up to translations. Using this, we can characterize the $\chi$-semistable locus of $T^*\mathbb{C}^n$ as follows.
\begin{proposition}[prop:ss_locus]{}
    A point $ (z, w) \in T^*\mathbb{C}^n $ is $\chi$-semistable if and only if
    \begin{equation*}
        \bigcap_{z_i = 0} H_i^+ \cap \bigcap_{w_i = 0} H_i^- \neq \emptyset,
    \end{equation*}
    where $ H_i^+ = \{ x \in \mathfrak{t}_{\mathbb{R}}^* \mid \pi^*(x)_i \geq \theta_i \} $ and $ H_i^- = \{ x \in \mathfrak{t}_{\mathbb{R}}^* \mid \pi^*(x)_i \leq \theta_i \} $.
\end{proposition}

A hyperplane arrangement is said to be
\begin{itemize}
    \item \vocab{simple} if every $k$ hyperplanes with non-empty intersection intersect at codimension $k$;
    \item \vocab{unimodular} if every collection of $d$ independent vectors in $\{\pi(e_1),\dots ,\pi(e_n)\}$ spans $\mathfrak{t}_{\mathbb{Z}}$, where $e_i$ denotes the standard generator of $\mathfrak{d}_{\mathbb{Z}}$; and
    \item \vocab{smooth} if it is simple and unimodular.
\end{itemize}
\inquiry[o]{I followed Mcbreen's quantum cohomology's notation}
The hyperplane arrangement can tell us how singular the hypertoric variety is.
\begin{proposition}[prop:smth_ht]{\cite[Theorem 3.2]{bielawski2000geometry}}
    The hypertoric variety $\mathcal{M}_{\chi, 0}$ has at worse orbifold singularities if and only if $\mathcal{A}_\chi$ is simple, and is smooth if and only if $\mathcal{A}_\chi$ is smooth.
\end{proposition}

Note that whether $\mathcal{A}_\chi$ is unimodular only depends on the exact sequence of tori, while whether it is simple depends on the character $\chi$. In fact, we can characterize this algebraically.
\begin{lemma}[lem:simple_iff_generic]{}
    The hyperplane arrangement $\mathcal{A}_\chi$ is simple if and only if for all $ x \in (\iota^*)^{-1}(\chi) $, $\{ \pi(e_i) \mid 1 \leq i \leq n, x_i = 0 \}$ is linearly independent in $\mathfrak{t}$, where $e_i$ denotes the $i$-th coordinate vector in $\mathfrak{d}$.
\end{lemma}

We say that $\chi$ is generic if the above condition is satisfied. 
In this case, the $G$ action on $\mu^{-1}(0)$ is locally free, and the semistable locus coincides with the stable locus. Any such choice of character gives rise to $T$-equivariantly diffeomorphic hypertoric varieties $\mathcal{M}_{\chi, 0}$ (see \cite[Lemma 2.1]{harada2004properties} and \cite[Theorem 6.1]{konno2003variation}).

Recall that for a graph $\Gamma$, we have a differential $ d_{\Gamma, \mathbb{C}^\times} : (\mathbb{C}^\times)^V \to (\mathbb{C}^\times)^E $, giving the following short exact sequence.
\begin{equation*}
    \begin{tikzcd}
        1 \arrow[r] & \im(d_{\Gamma, \mathbb{C}^\times}) \arrow[r, "\iota_\Gamma"] & (\mathbb{C}^\times)^E \arrow[r, "\pi_\Gamma"] & \HH^1(\Gamma, \mathbb{C}^\times) \arrow[r] & 1
    \end{tikzcd}
\end{equation*}
Denote by $G_\Gamma$, $D_\Gamma$ and $T_\Gamma$ the complex tori above, and by $ \mu_\Gamma : T^*\mathbb{C}^E \to \mathfrak{g}_\Gamma^* $ the moment map. For any choice of $ \chi \in \mathfrak{g}_{\Gamma, \mathbb{Z}}^* $ and $ \lambda \in \mathfrak{g}_\Gamma^* $, we denote by $\Hfinchar{\chi}$ the hyperplane arrangement in $\mathfrak{t}_{\Gamma, \mathbb{R}}^*$ and $ \Mfinchar{\chi, \lambda} = T^*\mathbb{C}^E \HQ[\chi, \lambda] G_\Gamma = \mu_\Gamma^{-1}(\lambda) \GIT[\chi] G_\Gamma $ the hypertoric variety associated to the above short exact sequence.

Note that $\pi_\Gamma$ is simply the pairing defined in \Cref{sec:graph_homology_and_cohomology}, which is unimodular. The associated hyperplane arrangement is also unimodular.
Since $\pi_\Gamma$ is a dual representation of $\Gamma$, \Cref{lem:simple_iff_generic} then translates to that $\Mfinchar{\chi, 0}$ is smooth if and only if for all $ x \in (\iota_{\mathbb{R}}^*)^{-1}(\chi) $, $ \{ e \in E \mid x_e = 0 \} \in \face$.
For any such choice of $\chi$, 
we will simply denote the smooth hypertoric variety by $\Mfin$. Also, $ \pi_\Gamma^* : \mathfrak{t}_{\Gamma}^* \to \mathfrak{d}_{\Gamma}^* $ is just the inclusion $ \HH_1(\Gamma, \mathbb{C}) \to \mathbb{C}^E $. Thus, $ \mu_\Gamma^{-1}(0) = \{ (z, w) \in T^*\mathbb{C}^E \mid (z_ew_e)_{e \in E} \in \HH_1(\Gamma, \mathbb{C}) \} $.

\begin{example}[exp:]{}

    Consider the graph in \Cref{fig:graph_TP2} with a choice of $\theta$ such that the associated hyperplane arrangement is \Cref{fig:hyperplane_TP2}.
    By \Cref{prop:ss_locus}, the unstable locus is given by $\{z_x=z_y=z_z=0\} \subseteq T^*\mathbb{C}^E$.
    The associated hypertoric variety is isomorphic to $T^*\mathbb{P}^2$.
\end{example}


Notice that we can generalize the above construction to matroids and get a hypertoric variety. The map $\pi$ will then be a representation of the matroid. By construction, the matroid can be written as a direct sum of two matroids if and only if the short exact sequence of complex tori can be written as a direct product of two. In this case, the associated hypertoric variety is simply the direct product of the two. 
In particular, if $\Gamma=\Gamma_1\wedge \Gamma_2$, i.e. $\Gamma$ is obtained by gluing $\Gamma_1$ and $\Gamma_2$ at a vertex, then $\Mfin\cong\Mfin[\Gamma_1]\times \Mfin[\Gamma_2]$.



\subsection{Cohomology of Hypertoric Varieties}

One can associate with any simplicial complex a commutative ring, called the Stanley-Reisner ring. We apply this construction to the simplicial coindependence complex $\face$ (cref. \Cref{def:face_coindependent_complex}) for any graph $\Gamma$.
\begin{definition}[def:SRfin]{}
    Let $\ZZ[E]$ be the polynomial ring generated by the edges of the graph $\Gamma$, where each generator $ e \in E $ is of degree $2$. Let $\idealI$ be the monomial ideal of $\ZZ[E]$ generated by bonds, i.e.
    \begin{equation*}
        \idealI = \left\langle \prod_{e \in S} e \mid S \subseteq E \text{ is a bond} \right\rangle \subseteq \ZZ[E].
    \end{equation*}
    Define the Stanley-Reisner ring associated with the graph $\Gamma$ by $ \SR = \ZZ[E] / \idealI $, which inherits the grading.
\end{definition}

If we further quotient out the ideal generated by $\HH_1(\Gamma)$, we obtain the following ring.
\begin{definition}[def:Rfin]{}
    Define $\Rfin$ to be the quotient of $\SR$ by the ideal generated by the image of the composite $ \HH_1(\Gamma) \hookrightarrow \ZZ^E \hookrightarrow \ZZ[E] \twoheadrightarrow \SR $. In other words, we have the following exact sequence.
    \begin{equation*}
        \begin{tikzcd}
            \SR[\bullet-2] \otimes \HH_1(\Gamma) \arrow[r] & \SR \arrow[r] & \Rfin \arrow[r] & 0
        \end{tikzcd}
    \end{equation*}
\end{definition}

\begin{theorem}[thm:eq_coho_ha]{\cite[Theorem 3.1]{konno2000cohomology}}
    For any generic character $\chi$, the $T_\Gamma$-equivariant cohomology of $\Mfinchar{\chi, 0}$ is isomorphic to $\SR$, and the ordinary cohomology is isomorphic to $\Rfin$.
\end{theorem}

The isomorphism could be understood via the following diagram, where the labels of the arrows denote the kernels of the maps.
\begin{equation} \label{eqn:cohom_ht}
    \begin{tikzcd}
        \ZZ[E] \arrow[r, phantom, "\cong"] & \HH_{D_\Gamma}^\bullet(T^*\mathbb{C}^E) \arrow[r, "{\idealI}"] \arrow[d, "{\HH_1(\Gamma)}"] & \HH_{D_\Gamma}^\bullet(\mu_\Gamma^{-1}(0)^{\chi\text{-ss}}) \arrow[r, "\sim"] \arrow[d, "{\HH_1(\Gamma)}"] & \HH_{T_\Gamma}^\bullet(\Mfinchar{\chi, 0}) \arrow[d, "{\HH_1(\Gamma)}"] \\
        {} & \HH_{G_\Gamma}^\bullet(T^*\mathbb{C}^E) \arrow[r, "{\idealI}"'] & \HH_{G_\Gamma}^\bullet(\mu_\Gamma^{-1}(0)^{\chi\text{-ss}}) \arrow[r, "\sim"'] & \HH^\bullet(\Mfinchar{\chi, 0})
    \end{tikzcd}
\end{equation}

\subsection{A Cochain Complex \texorpdfstring{$\HTfin$}{HTfin}}

We now introduce a cochain complex whose cohomology is isomorphic to $\Rfin$.

\begin{definition}[def:HTfin]{}
    Define a complex $(\HTfin, \dHTfin)$ by
    \begin{equation*}
        \HTfin[2p][q] = \bigoplus_{S \in \face[p]} \Ext[q] \HH_1(\Gamma \del S).
    \end{equation*}
    Given $ S \in \face[p] $ and $ \gamma_1, \dots, \gamma_q \in \HH_1(\Gamma \del S) $, the corresponding element in $\HTfin[2p+q]$ is denoted by $1_S \gamma_1 \dots \gamma_q$. The differential $ \dHTfin : \HTfin \to \HTfin[\bullet+2][\bullet-1] $ restricted to the summand $\Ext \HH_1(\Gamma \del S) \to \Ext[\bullet-1] \HH_1(\Gamma \del S \del e)$ is the interior product $\intp{e}$ by $ [e] \in \HH^1(\Gamma \del S) $.

    Define $ \grHTfin[k][p] = \HTfin[2p][k-p] $. Then, each $\grHTfin[k]$ is a subcomplex with differential of degree $1$, giving a decomposition $ \HTfin \cong \bigoplus_{k \in \mathbb{N}} \grHTfin[k]{[-k]} $.
\end{definition}

Given a coherent cotree $\cotreefunction$, we identify $ \HH_1(\Gamma \del S) \cong \ZZ^{\cotree{S}} \cong \HH^1(\Gamma \del S) $ for $ S \in \face[p] $. Given $ x_1, \dots, x_q \in \cotree{S} $ distinct, we have $ x_1 \dots x_q \in \Ext[k] \HH_1(\Gamma \del S) $, thus giving an element
\begin{equation*}
    1_S x_1 \dots x_q \in \HTfin[2p][q].
\end{equation*}
We will show algebraically that this complex computes the cohomology of $\Mfin$. The following theorem shows that we can use the subcomplex $\HTfin[\bullet][\leq 1]$ to recover $\Rfin$.

\begin{theorem}[thm:Rfin_pres]{}
    There is an exact sequence of abelian groups as follows.
    \begin{equation*}
        \begin{tikzcd}
            \HTfin[2k-2][1] = \bigoplus_{S \in \face[k-1]} \HH_1(\Gamma \del S) \arrow[r, "\dHTfin"] & \HTfin[2k][0] = \ZZ^{\face[k]} \arrow[r] & \Rfin[2k] \arrow[r] & 0
        \end{tikzcd}
    \end{equation*} 
\end{theorem}

The proof is obtained by combining the following two lemmas. We first establish some notations. Given $ S \in \face $, we write the corresponding generator in $\ZZ^{\face}$ as $1_S$. Given $ \sigma : E \to \NN $, write the corresponding monomial in $\ZZ[E]$ as $ E^\sigma = \prod_{e \in E} e^{\sigma(e)} $. Note that
\begin{equation*}
    \idealI = \langle E^\sigma \mid \supp(\sigma) \notin \face \rangle.
\end{equation*}
We say that $E^\sigma$ is square-free if $ \im \sigma \subseteq \{ 0, 1 \} $. There is an inclusion $ \ZZ^{\face} \hookrightarrow \SR $ as square-free monomials.

\begin{lemma}[lem:surj_face]{}
    The map $ \ZZ^{\face} \to \Rfin $ is surjective. In other words, $\Rfin$ is generated by square-free monomials.
    \begin{proof}
        We borrow the idea in the proof of Lemma A.4 in \cite{mcbreen2014intersection}. Pick a linear order on $E$. For all $E^\sigma$ that is not square-free, i.e. $ \sigma(e) > 1 $ for some $ e \in E $, define $ P(E^\sigma) = \max \{ e \in E \mid \sigma(e) > 1 \} $ and $ Q(E^\sigma) = \sigma(P(E^\sigma)) $. Define a preorder on the set of all monomials $E^\sigma$ so that $ E^\sigma \leq E^{\sigma'} $ is true if
        \begin{itemize}
            \item $E^\sigma$ is square-free, or
            \item $E^{\sigma'}$ is not square-free and $ P(E^\sigma) < P(E^{\sigma'}) $, or
            \item $E^{\sigma'}$ is not square-free, $ P(E^\sigma) = P(E^{\sigma'}) $ and $ Q(E^\sigma) \leq Q(E^{\sigma'}) $.
        \end{itemize}
        Since $\Rfin$ is a quotient of $\SR$, the former is generated by the monomials $E^\sigma$ with $ \supp(\sigma) \in \face $. We perform an induction, with respect to the preorder defined above, to show that every such monomial is in the image of $ \ZZ^{\face} \to \Rfin $. If $E^\sigma$ is minimal with respect to the preorder, i.e. it is square-free, then $E^\sigma$ is the image of $ 1_{\supp(\sigma)} \in \ZZ^{\face} $. Suppose $E^\sigma$ is not square-free. Write $ S = \supp(\sigma) $, $ e = P(E^\sigma) $ and $ E^\sigma = E^{\sigma'} \cdot e $. Since $e$ is not a bond in $\Gamma\del(S\setminus e)$, there exists $ \gamma_e \in \HH_1(\Gamma \del (S \setminus e)) $ such that $ \pair{\gamma_e}{e}[\Gamma \del (S \setminus e)] = 1 $. Since $ E^{\sigma'} \cdot \gamma_e $ is in the kernel of $ \SR \to \Rfin $,
        \begin{equation*}
            E^\sigma = E^{\sigma'} \cdot e = E^{\sigma'} \cdot (e - \gamma_e) = -\sum_{e' \in E \setminus S} \pair{\gamma_e}{e'}[\Gamma \del (S \setminus e)] E^{\sigma'} \cdot e'
        \end{equation*}
        in $\Rfin$. Note that for all $ e' \in E \setminus S $, since $ e' \notin \supp(\sigma') $, $ P(E^{\sigma'} \cdot e') = P(E^{\sigma'}) $. Thus, we either have
        \begin{itemize}
          \item $E^{\sigma'}$ is square-free, or
            \item $ P(E^{\sigma'}) < P(E^\sigma) $ if $ Q(E^\sigma) = 2 $, or
            \item $ P(E^{\sigma'}) = P(E^\sigma) $ and $ Q(E^{\sigma'}) = Q(E^\sigma) - 1 $ if $ Q(E^\sigma) > 2 $.
        \end{itemize}
        In any case, $ E^{\sigma'} \cdot e' < E^\sigma $. The induction is complete.
    \end{proof}
\end{lemma}

\begin{lemma}[lem:ker_face]{}
    The kernel of $ \ZZ^{\face} \to \Rfin $ equals the image of $ \dHTfin : \bigoplus_{S \in \face} \HH_1(\Gamma \del S) \to \ZZ^{\face} $.
    \begin{proof}
        Denote the map $ \ZZ^{\face} \to \Rfin $ by $\pi$ and the map $ \SR \to \Rfin $ by $\tau$. Compare the following short exact sequences.
        \begin{equation*}
            \begin{tikzcd}
                0 \arrow[r] & \ker \pi \arrow[r] \arrow[d, hook] & \ZZ^{\face} \arrow[r, "\pi"] \arrow[d, hook] & \Rfin \arrow[r] \arrow[d, equal] & 0 \\
                0 \arrow[r] & \ker \tau \arrow[r] & \SR \arrow[r, "\tau"] & \Rfin \arrow[r] & 0
            \end{tikzcd}
        \end{equation*}
        It follows that the left square is a pullback, i.e. $ \ker \pi = \ZZ^{\face} \cap \ker \tau $, where we identified $\ZZ^{\face}$ with its image in $\SR$. It is clear that $ \im \dHTfin \subseteq \ZZ^{\face} \cap \ker \tau $. It remains to show that $ \ZZ^{\face} \cap \ker \tau \subseteq \im \dHTfin $.
        
        We will use the setting in the previous lemma. Let $ \sum_{i=1}^k E^{\sigma_i} \cdot \gamma_i \in \ZZ^{\face} \cap \ker \tau $ so that the $E^{\sigma_i}$ are pairwise distinct. Consider the case where all the $E^{\sigma_i}$ are minimal, i.e. square-free. Suppose there exists $i$ such that $ E^{\sigma_i} \cdot \gamma_i \notin \im \dHTfin $, i.e. $ \supp(\sigma_i) \cap \supp(\gamma_i) \neq \emptyset $. Then, $ E^{\sigma_i} \cdot \gamma_i $ contains $ E^{\sigma_i} \cdot e $ as a summand. Since the whole sum is in $\ZZ^{\face}$, there exists $ j \neq i $ such that $ E^{\sigma_j} \cdot \gamma_j $ also contains $ E^{\sigma_i} \cdot e $ as a summand. Since $E^{\sigma_j}$ is square-free, it follows that $ E^{\sigma_i} = E^{\sigma_j} $, which contradicts our assumption. Thus, $ E^{\sigma_i} \cdot \gamma_i \in \im \dHTfin $ for all $i$.

        Now, consider the case where not all $E^{\sigma_i}$ are square-free, and let $E^{\sigma_i}$ be maximal among all the summands in the above preorder. Let $ e = P(E^{\sigma_i}) $. We claim that $ e' < e $ for all $ e' \in \supp(\sigma_i) \cap \supp(\gamma_i) $. Else, there will be $ j \neq i $ such that $ P(E^{\sigma_i}) = P(E^{\sigma_j}) $ and $ Q(E^{\sigma_i}) > Q(E^{\sigma_j}) $ (if $ e' = e $), or $ P(E^{\sigma_i}) > P(E^{\sigma_j}) $ (if $ e' > e $). In any case, $ E^{\sigma_i} > E^{\sigma_j} $, which contradicts maximality. Write $ E^{\sigma_i} = E^{\sigma_i'} \cdot e $ and pick $ \gamma_e \in \HH_1(\Gamma \del (S \setminus e)) $ such that $ \pair{\gamma_e}{e}[\Gamma \del (S \setminus e)] = 1 $. Then,
        \begin{align*}
            E^{\sigma_i} \cdot \gamma_i &= E^{\sigma_i'} \cdot e \cdot \gamma_i = E^{\sigma_i'} \cdot (e - \gamma_e) \cdot \gamma_i + E^{\sigma_i'} \cdot \gamma_i \cdot \gamma_e \\
            &= -\sum_{e' \in E \setminus S} \pair{\gamma_e}{e'}[\Gamma \del (S \setminus e)] (E^{\sigma_i'} \cdot e') \cdot \gamma_i + \sum_{e'' \in E} \pair{\gamma_i}{e''} (E^{\sigma_i'} \cdot e'') \cdot \gamma_e.
        \end{align*}
        For the first summand, we have shown in the previous proof that $ E^{\sigma_i'} \cdot e' < E^{\sigma_i} $. For the second summand, since $ e'' < e $ for all $ e'' \in \supp(\sigma_i) \cap \supp(\gamma_i) $, it follows that $ E^{\sigma_i'} \cdot e'' < E^{\sigma_i} $ when $ \pair{\gamma_i}{e''} \neq 0 $. Thus, this algorithm reduces the monomial preorder for all the maximal summands. By induction, we reduce any such sum to the square-free case, which the assertion then follows by the argument above.
    \end{proof}
\end{lemma}

\begin{example}[exp:face_surj_TP2]{}
    For the graph $\Gamma$ in \Cref{fig:graph_TP2}, we have
    \[
        \Rfin = \frac{\ZZ[x, y, z]}{xyz, x-z, y-z}.
    \]
    \Cref{lem:surj_face} says that any monomial in $\Rfin$ can be expressed as a sum of monomials corresponding to elements in $\face$.
    Note that the degree $1$ monomials $x$, $y$ and $z$ are all equal. Thus, any degree $2$ monomial equals to a square-free one, say $yz$. By a similar argument, any degree 3 monomial equals to $xyz$, which vanishes in $\Rfin$.
    Using the notations in the proof of \Cref{lem:ker_face}, we see that $z(x-z)-z(y-z) = zx-zy\in \ZZ^{\face} \cap \ker\tau$. \Cref{lem:ker_face} asserts that $zx-zy$ lies in $\im \dHTfin$. Indeed, $\HH_1(\Gamma\del z)$ is generated by $x-y=:\gamma_z$, hence $\dHTfin(1_z\gamma_z) = 1_{zx}-1_{zy}$.
\end{example}

\subsection{Showing the Homotopy Equivalence}

We will prove the following theorem in this subsection.

\begin{theorem}[thm:HTfin_eqv]{}
    There are homotopy equivalences
    \begin{equation*}
        \begin{tikzcd}
            (\grHTfin[k], \dHTfin) \arrow[r, "\fHTfin", shift left] & \arrow[l, "\gHTfin", shift left] \Rfin[2k]{[-k]}.
        \end{tikzcd}
    \end{equation*}
    \noskipline
    \begin{remark}
        In particular, the cohomology of $\HTfin$, graded by the total degree, is isomorphic to $\Rfin$ as graded abelian groups.
    \end{remark}
\end{theorem}

We will explicitly construct $\fHTfin$ and $\gHTfin$, as well as a homotopy $\hHTfin$ from the identity map to $\gHTfin\fHTfin$, depending on a choice of coherent cotree of $\Gamma$ and a further choice as indicated in \Cref{cor:J_basis}. We first fix a coherent cotree $\cotreefunction$ and show the following results.

\begin{lemma}[lem:cotree_elim]{}
    Let $ x, y \in \cotree{\emptyset} $. Then
    \begin{equation*}
        \dHTfin(1_x y) - \dHTfin(1_y x) = \sum_{t \in \tree{\emptyset}} \dHTfin (1_t (\pair{y}{t} x - \pair{x}{t} y)).
    \end{equation*}
    \noskipline
    \begin{proof}
        By \Cref{lem:tree_invariant}, $ \tree{x} = \tree{y} = \tree{\emptyset} $. In particular, $ x \in \cotree{y} $ and $ y \in \cotree{x} $. Using \Cref{eqn:hom_incl_in_coord}, both sides equate to
        \begin{equation*}
            \sum_{t \in \tree{\emptyset}} (\pair{y}{t} 1_{xt} - \pair{x}{t} 1_{yt}).
        \end{equation*}
    \end{proof}
\end{lemma}

\begin{corollary}[cor:cotree_elim]{}
    Let $ a_i \in \ZZ $, $ S_i \in \face $ and $ x_i \in \cotree{S_i} $ for all $ 1 \leq i \leq k $. Suppose $ \sum_{i=1}^k a_i = 0 $ and $ S_i \cup x_i = S_j \cup x_j $ for all $ 1 \leq i, j \leq k $. In particular, we have $ |\cotreein{S_i}| = |\cotreein{S_j}| $ for all $i$ and $j$ by \Cref{lem:cotreein_subset}. Then,
    \begin{equation*}
        \sum_{i=1}^k a_i \dHTfin(1_{S_i}x_i) = \sum_{i'=1}^{k'} a_{i'}' \dHTfin(1_{S_{i'}'}x_{i'}')
    \end{equation*}
    for some $a_{i'}'$, $S_{i'}'$ and $x_{i'}'$ such that $ |\cotreein{S_{i'}'}| < |\cotreein{S_i}| $ for all $i$ and $i'$.
    \begin{proof}
        For each $ i > 1 $, apply \Cref{lem:cotree_elim} to the graph $ \Gamma \del (S_1 \cap S_i) $ gives
        \begin{equation*}
            \dHTfin(1_{S_i}x_i) = \dHTfin(1_{S_1}x_1) + \sum_{t \in \tree{S_1 \cap S_i}} \left(\pair{x_i}{t}[\Gamma \del (S_1 \cap S_i)] \dHTfin(1_{(S_1 \cap S_i) \cup t} x_1) - \pair{x_1}{t}[\Gamma \del (S_1 \cap S_i)] \dHTfin(1_{(S_1 \cap S_i) \cup t} x_i)\right).
        \end{equation*}
        Moreover, we have $ \cotreein{(S_1 \cap S_i) \cup t} \subseteq \cotreein{S_1 \cap S_i} \subsetneq \cotreein{S_i} $ for $ t \in \tree{S_1 \cap S_i} $ by \Cref{lem:cotreein_subset}. The assumption $ \sum_{i=1}^k a_i = 0 $ guarantees that the term $\dHTfin(1_{S_1}x_1)$ will be canceled out after doing the above expansion for all $ i > 1 $.
    \end{proof}
\end{corollary}

We can then prove the following theorem.

\begin{theorem}[thm:Rfin_split]{}
    $ \ZZ^{\face} = \im \dHTfin \oplus \ZZ^{\basis} $, where $ \dHTfin : \bigoplus_{S \in \face} \HH_1(\Gamma \del S) \to \ZZ^{\face} $.
    \begin{remark}
        Equivalently, $\basis[k]$ gives a monomial basis for $\Rfin[2k]$ for all $ k \in \mathbb{N} $.
    \end{remark}
    \begin{proof}
        We first show that $\im \dHTfin$ and $\ZZ^{\basis}$ generate $\ZZ^{\face}$. Let $ S \in \face \setminus \basis $, i.e. $ \cotreein{S} \neq \emptyset $. Pick $ x \in \cotreein{S} $, i.e. $ x \in \cotree{S \setminus x} $. Then,
        \begin{equation*}
            1_S = \dHTfin(1_{S \setminus x}x) - \sum_{t \in \tree{S}} \pair{x}{t}[\Gamma \del (S \setminus x)] 1_{S \setminus x \cup t}.
        \end{equation*}
        Note that the first term is in $\im \dHTfin$ and for the second sum, we have $ \cotreein{S \setminus x \cup t} \subsetneq \cotreein{S} $ by \Cref{lem:cotreein_subset}. Thus, by induction, we can reduce all generators in $\ZZ^{\face}$ into a sum of elements in $\im \dHTfin$ and $\ZZ^{\basis}$.

        We now show that $ \im \dHTfin \cap \ZZ^{\basis} = \{ 0 \} $. Let
        \begin{equation*}
            \sum_{i=1}^k a_i \, \dHTfin(1_{S_i}x_i) = \sum_{i=1}^k a_i \left( 1_{S_i \cup x_i} + \sum_{t \in \tree{S_i}} \pair{x_i}{t}[\Gamma \del S_i] 1_{S_i \cup t} \right) \in \im \dHTfin \cap \ZZ^{\basis}
        \end{equation*}
        where $ a_i \in \ZZ $, $ S_i \in \face $ and $ x_i \in \cotree{S_i} $, so that the $1_{S_i}x_i$ are pairwise distinct as elements in $\bigoplus_{S \in \face} \HH_1(\Gamma \del S)$. Pick $i$ such that $\cotreein{S_i \cup x_i}$ is maximal. Since $x_i\in \cotreein{S_i\cup x_i}$, if $a_i\neq 0$, $1_{S_i\cup x_i}$ must be a summand of $\dHTfin(1_{S_j}x_j)$ for some $j\neq i$ for the final sum to be in $\mathbb{Z}^{\basis}$. Also, if $1_{S_i \cup x_i}$ is a summand of $\dHTfin(1_{S_j}x_j)$, we must have $ S_i \cup x_i = S_j \cup x_j $. Else, we will have $ \cotreein{S_i \cup x_i} \subsetneq \cotreein{S_j \cup x_j} $ violating maximality.

        Consider the base case where $ S_i \in \basis $ for all $i$, $ S_i \cup x_i = S_j \cup x_j $ would implies $ S_i = S_j $ and $ x_i = x_j $ as $x_i\in \cotreein{S_i\cup x_i\del e}$, i.e. $S_i\cup x_i\del e\notin \basis$, for all $e\in S_j$. Hence $\sum_{i=1}^{k}a_i\dHTfin(1_{S_i}x_i)\in \im \dHTfin\cap \mathbb{Z}^{\basis}$ if and only if $a_i=0$ for all $i$.

        Suppose $ S_i \notin \basis $ for some $i$. Then, we can partition $ \{ i \mid \cotreein{S_i \cup x_i} \text{ is maximal} \} $ according to the relation $ S_i \cup x_i = S_j \cup x_j $. Since the whole sum is in $\ZZ^{\basis}$, the $1_{S_i \cup x_i}$ must all be cancelled out. Thus, $ \sum a_i = 0 $ where $i$ ranges over any such partition. Applying \Cref{cor:cotree_elim} to all such partitions, by induction, the statement reduces to the base case where $ S_i \in \basis $ for all $i$. This concludes the proof.
        \inquiry[a]{Actually what does partition mean here? Not all $S_i\cup x_i$ is maximal here, right? So you mean you can pick a subset out?}
    \end{proof}
\end{theorem}

\begin{corollary}[cor:J_basis]{}
    Any choices of $ \choice{S} \in \cotreein{S} $ for each $ S \in \face \setminus \basis $ gives a basis $ \{ \dHTfin(1_Sx) \mid \choice{S \cup x} = x \} $ of $\im \dHTfin$.
    \begin{proof}
        Modifying the first part of the proof of \Cref{thm:Rfin_split} so that we pick $ \choice{S} \in \cotreein{S} $ for all $ S \in \face \setminus \basis $, we can see that $ \langle \dHTfin(1_Sx) \mid \choice{S \cup x} = x \rangle $ generates $\im \dHTfin$. Moreover, $ \dim \im \dHTfin = \dim \ZZ^{\face} - \dim \ZZ^{\basis} = |\face \setminus \basis| $, which is the cardinality of the asserted basis.
    \end{proof}
\end{corollary}

Extend the map $ \ZZ^{\face} \to \Rfin $ in \Cref{thm:Rfin_pres} into $ \fHTfin : \HTfin \to \Rfin $, and extend the split constructed in \Cref{thm:Rfin_split} into $ \gHTfin : \Rfin \to \HTfin $ so that they are trivial with respect to $\HTfin[\bullet][{>0}]$. It follows that $ \fHTfin\gHTfin = \mathrm{id} $. It is clear that they are chain maps with respect to the differential $\dHTfin$ on $\HTfin$ and the zero differential on $\Rfin$.

For each $ S \in \face \setminus \basis $, choose $ \choice{S} \in \cotreein{S} $ as in \Cref{cor:J_basis}. Define the homotopy $ \hHTfin : \HTfin \to \HTfin[\bullet-2][\bullet+1] $ recursively by
\begin{equation*}
    1_S x_1 \dots x_k \mapsto
    \begin{cases}
        0 &\quad \text{if } k=0 \text{ and } S \in \basis, \\
        0 &\quad \text{if } \choice{S \cup \{x_1, \dots, x_k\}} \notin S, \text{ and} \\
        1_{S \setminus x_0} x_0 \dots x_k - \sum_{t \in \tree{S \setminus x_0}} \hHTfin(1_{S \setminus x_0 \cup t} \intp{t}(x_0 \dots x_k)) &\quad \text{if } x_0 = \choice{S \cup \{x_1, \dots, x_k\}} \in S.
    \end{cases}
\end{equation*}
The recursion is well-defined since both $ \cotreein{S\del x_0},\cotreein{S \setminus x_0 \cup t} \subsetneq \cotreein{S} $ by \Cref{lem:cotreein_subset}.

We now show that $ \mathrm{id} - \gHTfin\fHTfin = \hHTfin\dHTfin + \dHTfin\hHTfin $.

\begin{lemma}[lem:]{}
    $(\mathrm{id} - \hHTfin\dHTfin - \dHTfin\hHTfin)(1_S x_1 \dots x_k) = 0 = \gHTfin\fHTfin(1_S x_1 \dots x_k)$ when $ k > 0 $.
    \begin{proof}
        We compute $\hHTfin\dHTfin(1_S x_1 \dots x_k)$ and $\dHTfin\hHTfin(1_S x_1 \dots x_k)$ by considering two cases.

        The first case is when $ x_r = \choice{S \cup \{x_1, \dots, x_k\}} \notin S $ for some $ 1 \leq r \leq k $. Then, $ \hHTfin(1_S x_1 \dots x_k) = 0 $. Also,
        \begin{align*}
            \hHTfin \dHTfin (1_S x_1 \dots x_k) &= \sum_{i=1}^k (-1)^{i-1} \hHTfin(1_{S \cup x_i} x_1 \dots \widehat{x_i} \dots x_k) + \sum_{t \in \tree{S}} \hHTfin(1_{S \cup t} \intp{t}(x_1 \dots x_k)) \\
            &= 1_S x_1 \dots x_k - \sum_{t \in \tree{S}} \hHTfin(1_{S \cup t} \intp{t}(x_1 \dots x_k)) + \sum_{t \in \tree{S}} \hHTfin(1_{S \cup t} \intp{t}(x_1 \dots x_k)) \\
            &= 1_S x_1 \dots x_k.
        \end{align*}

        Now we turn to the second case where $ x_0 = \choice{S \cup \{x_1, \dots, x_k\}} \in S $. Using \Cref{lem:tree_invariant}, we can show that
        \begin{equation*}
            \tree{S \setminus x_0 \cup x_i} = \tree{S \setminus x_0} = \tree{S}.
        \end{equation*}
        Thus,
        \begin{align*}
            \hHTfin \dHTfin (1_S x_1 \dots x_k) &= \sum_{i=1}^k (-1)^{i-1} \hHTfin(1_{S \cup x_i} x_1 \dots \widehat{x_i} \dots x_k) + \sum_{t \in \tree{S}} \hHTfin(1_{S \cup t} \intp{t}(x_1 \dots x_k)) \\
            &= \sum_{i=1}^k (-1)^{i-1} \left( 1_{S \setminus x_0 \cup x_i} x_0 \dots \widehat{x_i} \dots x_k - \sum_{t \in \tree{S \setminus x_0 \cup x_i}} \hHTfin(1_{S \setminus x_0 \cup x_i \cup t} \intp{t}(x_0 \dots \widehat{x_i} \dots x_k)) \right) \\
            &\quad + \sum_{t \in \tree{S}} \hHTfin(1_{S \cup t} \intp{t}(x_1 \dots x_k)) \\
            &= \sum_{i=1}^k (-1)^{i-1} 1_{S \setminus x_0 \cup x_i} x_0 \dots \widehat{x_i} \dots x_k + \sum_{i=0}^k (-1)^i \sum_{t \in \tree{S}} \hHTfin(1_{S \setminus x_0 \cup x_i \cup t} \intp{t}(x_0 \dots \widehat{x_i} \dots x_k));
        \end{align*}
        \begin{align*}
            \dHTfin \hHTfin (1_S x_1 \dots x_k) &= \dHTfin(1_{S \setminus x_0} x_0 \dots x_k) - \sum_{t \in \tree{S \setminus x_0}} \dHTfin\hHTfin(1_{S \setminus x_0 \cup t} \intp{t}(x_0 \dots x_k)) \\
            &= \sum_{i=0}^k (-1)^i 1_{S \setminus x_0 \cup x_i} x_0 \dots \widehat{x_i} \dots x_k + \sum_{t \in \tree{S \setminus x_0}} 1_{S \setminus x_0 \cup t} \intp{t}(x_0 \dots x_k) \\
            &\quad - \sum_{t \in \tree{S \setminus x_0}} \dHTfin\hHTfin(1_{S \setminus x_0 \cup t} \intp{t}(x_0 \dots x_k)) \\
            &= \sum_{i=0}^k (-1)^i 1_{S \setminus x_0 \cup x_i} x_0 \dots \widehat{x_i} \dots x_k + \sum_{t \in \tree{S}} (\mathrm{id} - \dHTfin\hHTfin) (1_{S \setminus x_0 \cup t} \intp{t}(x_0 \dots x_k)).
        \end{align*}
        Note that
        \begin{align*}
            0 &= \dHTfin^2(1_{S \setminus x_0} x_0 \dots x_k) = \sum_{i=0}^k (-1)^i \dHTfin(1_{S \setminus x_0 \cup x_i} x_0 \dots \widehat{x_i} \dots x_k) + \sum_{t \in \tree{S \setminus x_0}} \dHTfin(1_{S \setminus x \cup t} \intp{t}(x_0 \dots x_k)) \\
            &= \sum_{i=0}^k (-1)^i \left( \sum_{j=0}^{i-1} (-1)^j 1_{S \setminus x_0 \cup x_i \cup x_j} x_0 \dots \widehat{x_j} \dots \widehat{x_i} \dots x_k + \sum_{j=i+1}^k (-1)^{j-1} 1_{S \setminus x_0 \cup x_i \cup x_j} x_0 \dots \widehat{x_i} \dots \widehat{x_j} \dots x_k \right) \\
            &\quad + \sum_{i=0}^k (-1)^i \sum_{t \in \tree{S \setminus x_0 \cup x_i}} 1_{S \setminus x_0 \cup x_i \cup t} \intp{t}(x_0 \dots \widehat{x_i} \dots x_k) + \sum_{t \in \tree{S \setminus x_0}} \dHTfin(1_{S \setminus x \cup t} \intp{t}(x_0 \dots x_k)) \\
            &= \sum_{i=0}^k (-1)^i \sum_{t \in \tree{S}} 1_{S \setminus x_0 \cup x_i \cup t} \intp{t}(x_0 \dots \widehat{x_i} \dots x_k) + \sum_{t \in \tree{S}} \dHTfin(1_{S \setminus x \cup t} \intp{t}(x_0 \dots x_k)).
        \end{align*}
        Using this, we see that
        \begin{equation*}
            (\mathrm{id} - \hHTfin \dHTfin - \dHTfin \hHTfin) (1_S x_1 \dots x_k) = - \sum_{t \in \tree{S}} (\mathrm{id} - \hHTfin \dHTfin - \dHTfin \hHTfin) (1_{S \setminus x_0 \cup t} \intp{t}(x_0 \dots x_k)).
        \end{equation*}
        By induction on the cardinality of $\cotreein{S}$, using \Cref{lem:cotreein_subset}, we can reduce to the base case where $ \cotreein{S} = \emptyset $. Then, $ \choice{S \cup \{x_1, \dots, x_k\}} \notin S $, so by the first case, the assertion is true.
    \end{proof}
\end{lemma}

\begin{lemma}[lem:]{}
    $(\mathrm{id} - \hHTfin\dHTfin - \dHTfin\hHTfin)(1_S) = \gHTfin\fHTfin(1_S)$.
    \begin{proof}
      Note that $ \dHTfin(1_S) = 0 $ by definition. If $ S \in \basis $, then $ \hHTfin(1_S) = 0 $. Suppose $ \cotreein{S} \neq \emptyset $ and let $ x = \choice{S} \in \cotreein{S} \subseteq S $. By \Cref{lem:tree_invariant}, $ \tree{S \setminus x} = \tree{S} $. Thus,
        \begin{align*}
            \dHTfin\hHTfin(1_S) &= \dHTfin(1_{S \setminus x} x) - \sum_{t \in \tree{S \setminus x}} \pair{x}{t}[\Gamma \del S] \dHTfin\hHTfin(1_{S \setminus x \cup t}) \\
            &= 1_S + \sum_{t \in \tree{S}} \pair{x}{t}[\Gamma \del S] (\mathrm{id} - \hHTfin \dHTfin - \dHTfin \hHTfin)(1_{S \setminus x \cup t}).
        \end{align*}
        We obtained the following recursive formula for $\mathrm{id} - \hHTfin \dHTfin - \dHTfin \hHTfin$ restricted to $ \ZZ^{\face} \to \ZZ^{\face} $.
        \begin{align*}
            1_S &\mapsto
            \begin{cases}
                1_S &\quad \text{if } S \in \basis, \text{ and} \\
                -\sum_{t \in \tree{S \setminus x}} \pair{x}{t}[\Gamma \setminus (S \setminus x)] (\mathrm{id} - \hHTfin \dHTfin - \dHTfin \hHTfin)(1_{S \setminus x \cup t}) &\quad \text{if } \choice{S} = x \in S.
            \end{cases}
        \end{align*}
        Note that the image of this map is $\ZZ^{\basis}$ since $ \cotreein{S \setminus x \cup t} \subsetneq \cotreein{S} $ by \Cref{lem:cotreein_subset} in each step of the recursion. It is clear that $\dHTfin(1_Sx)$ is in the kernel of this map if $ \choice{S \cup x} = x $. By \Cref{cor:J_basis}, the kernel of this map is contained in the image of $ \dHTfin : \bigoplus_{S \in \face} \HH_1(\Gamma \del S) \to \ZZ^{\face} $. It follows that $\mathrm{id} - \hHTfin \dHTfin - \dHTfin \hHTfin$ must be the projection $ \ZZ^{\face} \to \Rfin $ followed by the inclusion $ \Rfin \to \ZZ^{\face} $ along the split in \Cref{thm:Rfin_split}. Thus, $ \mathrm{id} - \hHTfin \dHTfin - \dHTfin \hHTfin = \gHTfin\fHTfin $.
        \inquiry[a]{Explain last sentence}
    \end{proof}
\end{lemma}

The above two lemmas, together with the description of $\fHTfin$ and $\gHTfin$, conclude \Cref{thm:HTfin_eqv}. Moreover, we have the following corollary.

\begin{corollary}[cor:Rfin_les]{}
    The homotopy equivalence gives a long exact sequence extending the one in \Cref{thm:Rfin_pres}.
    \begin{equation*}
        \begin{tikzcd}
            0 \arrow[r] & \grHTfin[k][0] \arrow[r, shift left, "\dHTfin"] & \cdots \arrow[l, shift left, dashed, "\hHTfin"] \arrow[r, shift left, "\dHTfin"] & \grHTfin[k][k-1] \arrow[l, shift left, dashed, "\hHTfin"] \arrow[r, shift left, "\dHTfin"] & \grHTfin[k][k] \arrow[l, shift left, dashed, "\hHTfin"] \arrow[r, shift left, "\fHTfin"] & \Rfin[2k] \arrow[r] \arrow[l, shift left, dashed, "\gHTfin"] & 0
        \end{tikzcd}
    \end{equation*}
\end{corollary}

\begin{example}{}
    Consider the graph $\Gamma$ in \Cref{fig:graph_TP2} with the coherent cotree $\cotreefunction$ defined in \Cref{exp:coherent_cotree_TP2}.
    The choice $\cotreefunction$ fixes a basis for $\Rfin$, namely $\{ 1, z, yz \}$. Hence
    $g: \Rfin \to \ZZ^{\face}$ is defined by sending all degree 0 monomials to $1_\emptyset$, all degree 1 monomials to $1_z$, all degree 2 monomials to $1_{yz}$, and all degree 3 or above monomials to 0.

    Fix the following choices $[S]\in \cotreein{S}$ for each $S\in\face\del\basis$ as in 
    \Cref{cor:J_basis}.
    \[
        \begin{array}{c|ccccccc}
            S\in \face & \emptyset & \CB{x} & \CB{y} & \CB{z} & \CB{x,y} & \CB{x,z} & \CB{y,z}\\
            \hline
            \cotreein{S} & \emptyset & \CB{x} & \CB{y} & \emptyset & \CB{x,y} & \CB{x} & \emptyset \\
            {[S]} &  & x & y &  & x & x & 
        \end{array}
    \]
    Write $\gamma_x = y-z$, $\gamma_y = x-z$ and $\gamma_z = x-y$, regarded as elements in $\HH_1(\Gamma)$.
    Then, $h:\HTfin \to \HTfin[\bullet -2][\bullet+1]$ is defined as follows.
    \[
        \begin{array}{c|cccccccccc}
            \alpha \in \HTfin & 1_\emptyset & 1_x & 1_y & 1_z & 1_{x}\gamma_x & 1_{y}\gamma_y & 1_{z}\gamma_z & 1_{xy} & 1_{xz} & 1_{yz} \\
            \hline
            h(\alpha) & 0 & \gamma_y & \gamma_x & 0 & \gamma_y\wedge\gamma_x & \gamma_x\wedge\gamma_y & 0 & 1_y\gamma_y & 1_z\gamma_z & 0
        \end{array}
    \]
\end{example}

\subsection{Poincaré Polynomial of Hypertoric Varieties}
\label{sec:poincare_polynomial_hypertoric}

The Tutte polynomial of a graph $\Gamma$ is defined by
\[
    \Tutte{x,y}=\sum_{S\subseteq E}^{}(x-1)^{k(S)-k(E)}(y-1)^{k(S)+|S|-|V|},
\]
where $k(S)$ denotes the number of connected components of the subgraph of $\Gamma$ with edge $S$.
An equivalent definition of it is
\[
    \Tutte{x,y}=\sum_{T\text{ spanning tree}}^{}x^{|\internalactivity{T}|}y^{|\externalactivity{T}|},
\]
where $\internalactivity{T}$ denotes internal activity of $T$, the dual concept of external activity. It is a non-trivial combinatorial fact that the second definition, which depends on the ordering on $E$, equals to the first definition that does not.
The $h$-polynomial of the dual graph of $\Gamma$ is a specialization of the Tutte polynomial defined by
\[
    \hpoly{q} := \Tutte{1,q} = 
    \sum_{T\text{ spanning tree}}^{}q^{|\externalactivity{T}|}.
\]

It is proven in \cite{hausel2002toric} and \cite{crawley2004absolutely} that the Poincaré polynomial of $\Mfin$ equals $\hpoly{q}$, i.e.
\begin{equation*}
    \hpoly{q} = \sum_{k=0}^d \dim_\ZZ \HH^{2(d-k)}(\Mfin, \ZZ) \, q^k = \sum_{k=0}^d \dim_\ZZ \Rfin[2(d-k)] \, q^k,
\end{equation*}
where $d$ is the genus of $\Gamma$.
As a corollary of \Cref{thm:Rfin_split} and \Cref{pps:coherent_cotree_and_external_activity}, if we choose a coherent cotree induced by an ordering on the set $E$, we obtain a monomial basis for $\Rfin$ labelled by spanning trees $T$ with degree $2(d - |\externalactivity{T}|)$. This basis agrees with the one given in \cite[Theorem III.2.5]{stanley1996combinatorics}.

On the other hand, \cite{abdelgadir2022tutte} also proved the above result by giving a stratification of hypertoric varieties
\begin{equation*}
    \Mfin=\coprod_{T\text{ spanning tree}}\mathbb{C}^{d+|\externalactivity{T}|}
\end{equation*}
depending on an ordering of $E$.

Recall from \Cref{pps:coherent_cotree_and_external_activity} that
if the coherent cotree is constructed from the same ordering on the set $E$, 
$\externalactivity{T} = \cotreein{T^* := E\del T}$.
We suspect that the Poincaré dual of the closure of the strata $\mathbb{C}^{d+|\externalactivity{T}|}$, which is of codimension $d - |\externalactivity{T}|$, is the cohomology class of $\Mfin$ of degree $2(d - |\externalactivity{T}|)$ given by the basis element $T^* \del \cotreein{T^*} \in \basis$ up to a sign (cref. \Cref{crl:basis_and_cotreein}).
\minormissing{mellit's paper}

\subsection{Deletion and Contraction of Hypertoric Varieties}

Hypertoric varieties behave nicely with respect to deletion and contraction of graphs.

We first discuss deletion. Let $ S \subseteq E $ be a set of edges that does not contain any bonds, i.e. $ S \in \face $. Then, we have $ G_{\Gamma \del S} \cong G_\Gamma $ canonically.
\begin{lemma}[lem:ht_del]{}
    For any character $ \chi \in \mathfrak{g}_{\Gamma, \ZZ}^* \cong \mathfrak{g}_{\Gamma \del S, \ZZ}^* $, there is a closed embedding $ \Mfinchar{\chi, 0}[\Gamma \del S] \hookrightarrow \Mfinchar{\chi, 0} $.
\end{lemma}

The embedding is given by the coordinate inclusion $ T^*\mathbb{C}^{E \del S} \hookrightarrow T^*\mathbb{C}^E $ in the prequotient. This is clearly equivariant with respect to the action of $ G_{\Gamma \del S} \cong G_\Gamma $, and restricts to a closed embedding $ \mu_{\Gamma \del S}^{-1}(0)^{\chi\text{-ss}} \to \mu_{\Gamma}^{-1}(0)^{\chi\text{-ss}} $.

Contraction is more complicated. Let $ S \subseteq E $ be a set of edges that does not contain any cycle. Recall that for $ (z, w) \in T^*\mathbb{C}^{E \con S} $, $ \mu_{\Gamma \con S}(z, w) = 0 $ if and only if $ (z_ew_e)_{e \in E \con S} \in \HH_1(\Gamma \con S, \mathbb{C}) $. Since $ \HH_1(\Gamma \con S, \mathbb{C}) \cong \HH_1(\Gamma, \mathbb{C}) $, there is a unique $ \eta \in \mathbb{C}^S $ making $ ((z_ew_e)_{e \in E \con S}, \eta) \in \HH_1(\Gamma, \mathbb{C}) $. For any partition $ S = S_+ \sqcup S_- $, this defines a map $ \mu_{\Gamma \con S}^{-1}(0) \to \mu_{\Gamma}^{-1}(0) $ given by
\begin{equation*}
    (z_e, w_e)_{e \in E \con S} \mapsto ((z_e, w_e)_{e \in E \con S}, (1, \eta_e)_{e \in S_+}, (\eta_e, 1)_{e \in S_-}).
\end{equation*}
This map is equivariant with respect to the actions on both sides under the inclusion $ G_{\Gamma \con S} \hookrightarrow G_\Gamma $. Let $ \chi \in \mathfrak{g}_{\Gamma \con S, \ZZ}^* $. Since $ \mathfrak{t}_{\Gamma, \mathbb{R}}^* \cong \mathfrak{t}_{\Gamma \con S, \mathbb{R}}^* $, for any lift $\tilde{\chi}$ of $\chi$ along the projection $ \mathfrak{g}_{\Gamma, \ZZ}^* \to \mathfrak{g}_{\Gamma \con S, \ZZ}^* $, the ambient spaces of the hyperplane arrangements $\Hfinchar{\chi}$ and $\Hfinchar{\tilde{\chi}}[\Gamma \con S]$ are canonically identified. We impose a condition on $\tilde{\chi}$ that for all spanning cotree $ T^* \subseteq E \con S $ of $ \Gamma \con S $, the following holds in the hyperplane arrangement $\Hfinchar{\chi}$.
\begin{equation} \label{eq:cond_for_con}
    \bigcap_{e \in T^*} H_e \subseteq \bigcap_{e \in S_-} H_e^+ \cap \bigcap_{e \in S_+} H_e^-
\end{equation}
Note that $\bigcap_{e \in T^*} H_e$ is always a singleton. This ensures that the $\chi$-semistable points in $\mu_{\Gamma \con S}^{-1}(0)$ are mapped to $\tilde{\chi}$-semistable points in $\mu_{\Gamma}^{-1}(0)$. Such lift always exists since we have only finitely many conditions and we can move the hyperplanes $\{H_e\}_{e \in S}$ as far away as we want by varying $\tilde{\chi}$. Moreover, if $\chi$ is generic, we can always choose such lift to also be generic. Thus, we have the following result.
\begin{lemma}[lem:ht_con]{}
    For any character $ \chi \in \mathfrak{g}_{\Gamma \con S, \ZZ}^* $, there exists a lift $ \tilde{\chi} \in \mathfrak{g}_{\Gamma, \ZZ}^* $ with an open embedding $ \Mfinchar{\chi, 0}[\Gamma \con S] \hookrightarrow \Mfinchar{\tilde{\chi}, 0} $. Furthermore, we can choose $\tilde{\chi}$ to be generic if $\chi$ is.
\end{lemma}

To see that this is an open embedding, note that we can pick a spanning tree $ T \subseteq E $ containing $S$, giving independent generators $ g_e \in G_\Gamma $ for each $ e \in S $. For any $ (z, w) \in \mu_{\Gamma}^{-1}(0) $ with $ w_e \neq 0 $ for $ e \in S_+ $ and $ z_e \neq 0 $ for $ e \in S_- $, $ g \cdot (z, w) $ is in the image of the inclusion for $ g = \prod_{e \in S_+} w_eg_e \prod_{e \in S_-} z_e^{-1}g_e \in G_\Gamma $. Thus, the image of the inclusion is given by $ \{ [z, w] \in \Mfinchar{\tilde{\chi}, 0} \mid w_e \neq 0 \text{ for } e \in S_+ \text{ and } z_e \neq 0 \text{ for } e \in S_- \} $.
\minormissing{such lift always exists and generic?}

Assume now that $ e \in E $ is neither a bond nor a loop, $ S_+ = \emptyset $ and $ S_- = \{ e \} $. Denote the complement of the image of the contraction embedding by $D_e$, i.e. $ D_e = \{ [z, w] \in \Mfinchar{\tilde{\chi}, 0} \mid w_e = 0 \} $. In fact, the fundamental class of $D_e$ is the Poincaré dual of the cohomology class $ e \in \Rfin[2] \cong \HH^2(\Mfinchar{\tilde{\chi}, 0}) $. We claim that $D_e$ is a line bundle over the image of the deletion embedding. To see this, consider the coordinate projection $ T^*\mathbb{C}^E \to T^*\mathbb{C}^{E \del e} $, which clearly decends to $ \mu_\Gamma^{-1}(0) \cap \{ w_e = 0 \} \to \mu_{\Gamma \del e}^{-1}(0) $ since any cycle in $\Gamma$ that does not contain $e$ is a cycle in $\Gamma \del e$. Moreover, the fiber of any semistable point is semistable. This decends to a line bundle in the quotient. Thus, the embeddings $ \Mfinchar{\tilde{\chi}, 0}[\Gamma \del e] \hookrightarrow D_e \subseteq \Mfinchar{\tilde{\chi}, 0} \hookleftarrow \Mfinchar{\chi, 0}[\Gamma \con e] $ induces a Gysin sequence as follows.
\begin{equation*}
    \begin{tikzcd}
        \cdots \arrow[r] & \HH^{\bullet-2}(\Mfinchar{\tilde{\chi}, 0}[\Gamma \del e]) \arrow[r] & \HH^{\bullet}(\Mfinchar{\tilde{\chi}, 0}) \arrow[r] & \HH^{\bullet}(\Mfinchar{\chi, 0}[\Gamma \con e]) \arrow[r] & \cdots
    \end{tikzcd}
\end{equation*}

Since $\hpoly{q} = \hpoly[\Gamma \del e]{q} + \hpoly[\Gamma \con e]{q}$ is true by definition, we have $\dim_\ZZ \Rfin[2k] = \dim_\ZZ \Rfin[2k-2][\Gamma \del e] + \dim_\ZZ \Rfin[2k][\Gamma \con e]$. Thus, the connecting homomorphisms of the sequence above should be zero. From the point of view of \Cref{eqn:cohom_ht}, we should start from the following.
\begin{equation*}
    \begin{tikzcd}
        0 \arrow[r] & \ZZ[E \setminus e] \arrow[r] & \ZZ[E] \arrow[r] & \ZZ[E \con e] \arrow[r] & 0
    \end{tikzcd}
\end{equation*}
The inclusion is given by $ r \mapsto re $ and the projection is the quotient map realizing the isomorphism $ \ZZ[E \setminus e] \cong \ZZ[E] / \langle e \rangle $.

\begin{theorem}[thm:Rfin_delcon]{}
    If $ e \in E $ is neither a bond or a loop, there is a deletion-contraction short exact sequence.
    \begin{equation*}
        \begin{tikzcd}
            0 \arrow[r] & \Rfin[\bullet-2][\Gamma \del e] \arrow[r] & \Rfin \arrow[r] & \Rfin[\bullet][\Gamma \con e] \arrow[r] & 0
        \end{tikzcd}
    \end{equation*}
    \noskipline
    \begin{proof}
        We first show that both maps descend to the quotients.
        
        The map $\ZZ[E \setminus e] \to \ZZ[E]$ sends a monomial $(E \setminus e)^\sigma$ to the monomial $E^{\sigma'}$ so that $ \sigma'(e') = \sigma(e') $ if $ e' \in E \setminus e $ and $ \sigma'(e) = 1 $. Note that $ \supp(\sigma') = \supp(\sigma) \cup e $. Since $e$ is not a bond, $ \supp(\sigma) \in \face[\bullet][\Gamma \del e] $ if and only if $ \supp(\sigma) \cup e \in \face $. Also, the inclusion $ \ZZ^{E \setminus e} \hookrightarrow \ZZ^E $ restricts to an inclusion $ \HH_1(\Gamma \del e) \hookrightarrow \HH_1(\Gamma) $. It follows that we have a map $ \Rfin[\bullet-2][\Gamma \del e] \to \Rfin $.

        The map $\ZZ[E] \to \ZZ[E \setminus e]$ sends a monomial $E^\sigma$ to $(E \setminus e)^{\sigma|_{E \setminus e}}$ if $ \sigma(e) = 0 $, and $0$ else. Note that $ \supp(\sigma) \in \face[\bullet][\Gamma \con e] $ if and only if $ e \notin \supp(\sigma) $ and $ \supp(\sigma) \in \face $. Also, the projection $ \ZZ^E \twoheadrightarrow \ZZ^{E \setminus e} $ restricts to an isomorphism $ \HH_1(\Gamma) \cong \HH_1(\Gamma \con e) $ since $e$ is not a loop. It follows that we have a map $ \Rfin \to \Rfin[\bullet][\Gamma \con e] $.

        Finally, we show that the sequence is indeed short exact by expressing the maps in terms of bases. Recall that to any coherent cotree $\cotreefunction$ of $\Gamma$, we constructed a basis $\basis$ of $\Rfin$ in \Cref{thm:Rfin_split}. This also gives a basis
        \begin{equation*}
            \basis[\bullet][\Gamma \del e] = \{ S \subseteq E \setminus e \mid S \cup e \in \basis \}
        \end{equation*}
        of $\Rfin[\bullet][\Gamma \del e]$. Thus, the inclusion $ \Rfin[\bullet-2][\Gamma \del e] \to \Rfin $ is induced by the inclusion $ \basis[\bullet][\Gamma \del e] \hookrightarrow \basis $ of bases given by $ S \mapsto S \cup e $. Now assume that $ e \in \tree{\emptyset} $. Then, the coherent cotree $\cotreefunction[\Gamma \con e]$ gives a basis
        \begin{equation*}
            \basis[\bullet][\Gamma \con e] = \{ S \subseteq E \setminus e \mid S \in \basis \}.
        \end{equation*}
        Thus, the quotient $ \Rfin \to \Rfin[\bullet][\Gamma \con e] $ is given by the inclusion $ \basis[\bullet][\Gamma \con e] \subseteq \basis $ of bases. These inclusions give a partition of $\basis$, proving the exactness.
    \end{proof}
\end{theorem}

%% file: sec_periodization_new.tex
\section{Periodization}
\label{sec:periodizaiton}

\subsection{Periodization of Graphs}
\label{sec:periodization_of_graphs}

Write $ [a, b] = \{ i \in \ZZ \mid a \leq i \leq b \} $. Given a finite graph $ \Gamma = (V, E, h, t) $, we define a family of graphs $\{ \Gamma_n = (V_n, E_n, h_n, t_n) \}_{n \in \mathbb{N}}$ as follows. Let $ V_n = V \sqcup ([-n, n-1] \times E) $ and $ E_n = [-n, n] \times E $. Define $ h_n, t_n : E_n \to V_n $ by
\begin{equation*}
    h_n : (i, e) \mapsto
    \begin{cases}
        h(e) \quad &\text{if } i = n, \\
        (i, e) \quad &\text{else;}
    \end{cases}
    \quad \text{ and } \quad
    t_n : (i, e) \mapsto
    \begin{cases}
        t(e) \quad &\text{if } i = -n, \\
        (i-1, e) \quad &\text{else.}
    \end{cases}
\end{equation*}
In particular, $ \Gamma_0 \cong \Gamma $. Denote by $e_i$ the element $ (i, e) \in E_n $. Intuitively, $\Gamma_n$ is obtained by segmenting each edge $ e \in E $ into $2n+1$ edges $ e_{-n}, \dots, e_n \in E_n $. There is an evident contraction $ \Gamma_{n+1} \to \Gamma_n $ contracting the set of edges $ \{ e_{-n-1}, e_{n+1} \mid e \in E \} $. Heuristically, we would like to consider the inverse limit $ \Gamma_\per = \varprojlim_{n \in \mathbb{N}} \Gamma_n $, with $ E_\per = \ZZ \times E $ as the set of edges, as the periodization of $\Gamma$. The actual inverse limit in the category of reflexive graphs is not what we want. However, we could still consider $\Gamma_\per$ as an infinite matroid with evident infinite circuits, whose dual matroid is finitary. Practically, we will only work with the finite graphs $\Gamma_n$ to construct various other objects, and then take the appropriate limit afterward.

For $ S \in \face $ and $ I \in [-n, n]^S $, write $ S_I = \{ e_{I_e} \mid e \in S \} \subseteq E_n $. Since every pair of edges $e_i$ and $e_j$ with $i \neq j$ forms a bond, it is clear that these are all the cotrees. Thus, $ \face[\bullet][\Gamma_n] = \{ S_I \mid S \in \face, I \in [-n, n]^S \} $. It follows that the coindependent complex of $\Gamma_\per$ is given by the direct limit
\begin{equation*}
    \faceper = \varinjlim_{n \in \mathbb{N}} \face[\bullet][\Gamma_n] = \{ S_I \mid S \in \face, I \in \ZZ^S \}.
\end{equation*}

We construct a coherent cotree $\cotreefunction[\Gamma_n]$ out of $\cotreefunction$ by setting $ \cotree{S_I}[\Gamma_n] = \{ e_0 \mid e \in \cotree{S} \} $ for all $ S_I \in \face[\bullet][\Gamma_n] $. Coherence is immediate. To see that they are spanning cotrees, note that $\cotree{S_I}[\Gamma_n]$ does not contain any bond and $ |\cotree{S_I}[\Gamma_n]| = |\cotree{S}| = \dim \HH_1(\Gamma) = \dim \HH_1(\Gamma_n) $.

\begin{lemma}[lem:cotreein_periodization]{}
    For $\cotreefunction[\Gamma_n]$ as defined above and $ S_I \in \face[\bullet][\Gamma_n] $, we have
    \begin{equation*}
        \cotreein{S_I}[\Gamma_n] = \{ e_0 \mid e \in \cotreein{S} \cap (S \setminus \supp(I)) \}.
    \end{equation*}
    Consequently, $ \basis[\bullet][\Gamma_n] = \{ S_I \in \face[\bullet][\Gamma_n] \mid \cotreein{S} \subseteq \supp(I) \} $.
    \begin{proof}
        Write $ \supp(S_I) = \{ e_i \in S_I \mid i \neq 0 \} $. By the construction of $\cotreefunction[\Gamma_n]$, $ \supp(S_I) \subseteq \tree{S_I \setminus \supp(S_I)}[\Gamma_n] $. Thus, $ \cotreein{S_I}[\Gamma_n] = \cotreein{S_I \setminus \supp(S_I)}[\Gamma_n \del \supp(S_I)] $ by \Cref{lem:cotreein_subset}. Note that $ \{ e_i \mid e \in \supp(I), i \neq I_e \} \cup \{ e_i \mid e \in E \setminus \supp(I), i \neq 0 \} \subseteq \tree{\emptyset}[\Gamma_n \del \supp(S_I)] $. We can hence obtain a coherent cotree after contracting this set of edges. The desired form then follows from the identity $ \cotreein{S \setminus \supp(I)}[\Gamma \del \supp(I)] = \cotreein{S} \cap (S \setminus \supp(I)) $. 
    \end{proof}
\end{lemma}

\subsection{Periodized Hypertoric Variety}

Periodized hypertoric varieties are the periodic version of hypertoric varieties, which could be thought of as the hypertoric varieties associated with periodized hyperplane arrangements. Let $\Gamma$ be a finite graph, $ \chi \in \mathfrak{g}_{\Gamma, \ZZ}^* $ and $ \theta \in \mathfrak{d}_{\Gamma, \ZZ}^* $ be a lift of $\chi$ along $\iota_{\Gamma}^*$. Define $ \theta_n \in \mathfrak{d}_{\Gamma_n, \ZZ}^* $ by $ (\theta_n)_{e_i} = \theta_e + iN $ for some $ N \in \mathbb{N} $ and $ \chi_n = \iota_{\Gamma_n}^*(\theta_n) $. Assuming $\chi$ is generic, for large enough $N$, $\chi_n$ will be generic for all $ n \in \mathbb{N} $. This gives a family of hyperplane arrangements $\{\Hfinchar{\chi_n}[\Gamma_n]\}_{n \in \mathbb{N}}$, whose ambient spaces are all canonically identified with $\mathfrak{t}_{\Gamma, \ZZ}^*$. Taking the union, we obtain the periodized hyperplane arrangement, which is explicitly given by
\minormissing{why is $\chi_n$ generic}
\begin{equation*}
    \Hperchar{\chi} = \{ H_{e, i} = \{ x \in \mathfrak{t}_{\Gamma, \mathbb{R}}^* \mid \pair{x}{e} = \theta_e + iN \} \}_{e \in E, i \in \ZZ}.
\end{equation*}

For $ \Gamma = \loopgraph $ the loop graph, $ \chi = 0 $, $ \theta = 0 $ and $ N = 1 $, \Cref{lem:ht_con} gives us a sequence of open embeddings
\begin{equation*}
    \Mfinchar{\chi_0, 0}[\loopgraph_0] \hookrightarrow \Mfinchar{\chi_1, 0}[\loopgraph_1] \hookrightarrow \Mfinchar{\chi_2, 0}[\loopgraph_2] \hookrightarrow \cdots.
\end{equation*}
Denote by $\Mper[\loopgraph]$ the direct limit of the above sequence. Explicitly, it is given by the quotient of $ \ZZ \times T^*\mathbb{C} $ by the relation generated by $ (i, z, w) \sim (i+1, w^{-1}, zw^2) $ for all $ i \in \ZZ $, $ z \in \mathbb{C} $ and $ w \in \mathbb{C}^\times $. This space admits a $\mathbb{C}^\times$-action given by $ t \cdot [i, z, w] = [i, tz, t^{-1}w] $ and a moment map $ \mu^\per : \Mper[\loopgraph] \to \mathbb{C} $ given by $ [i, z, w] \mapsto zw $.

For a general graph $\Gamma$, we consider the subgroup $ G_\Gamma \subseteq D_\Gamma $ acting on $\Mper[\loopgraph]^E$ with a moment map $ \mu^\per_\Gamma : \Mper[\loopgraph]^E \to \mathfrak{d}_\Gamma^* \to \mathfrak{g}_\Gamma^* $. Together with $ \chi \in \mathfrak{g}_{\Gamma, \ZZ}^* $ and $ \lambda \in \mathfrak{g}_{\Gamma}^* $, we define the periodized hypertoric variety associated with these data to be the GIT quotient
\begin{equation*}
    \Mperchar{\chi, \lambda} = \Mper[\loopgraph]^E \HQ[\chi, \lambda] G_\Gamma = (\mu_\Gamma^\per)^{-1}(\lambda) \GIT[\chi] G_\Gamma.
\end{equation*}
We can characterize the $\chi$-semistable points of $\Mper[\loopgraph]^E$ using the periodized hyperplane arrangement $\Hperchar{\chi}$ as follows.
\minormissing{characterize ss points in Mper}
\begin{proposition}[]{}
    A point $ [i, z, w] \in \Mper[\loopgraph]^E $ is $\chi$-semistable if and only if
    \begin{equation*}
        \bigcap_{z'_e=0} H_{e, i'_e}^+ \cap \bigcap_{w'_e=0} H_{e, i'_e}^- \neq \emptyset
    \end{equation*}
    for all $ [i', z', w'] = [i, z, w] $.
\end{proposition}

\begin{example}[exp:periodized_loopgraph]{}
  For $\Gamma = \loopgraph$, the periodized hyperplane arrangement is precisely the lattice $\ZZ \subseteq \mathbb{R} \cong \mathfrak{t}^*_{\loopgraph, \mathbb{R}}$. The core of $\Mper[\loopgraph]$ is an infinite chain of $\mathbb{P}^1$, where consecutive $\mathbb{P}^1$'s are glued together along a single point, see \Cref{fig:inf_chain_P_one}.
\end{example}

\begin{minipage}{.3\linewidth}
  \begin{figure}[H]
    \centering
    \includegraphics[]{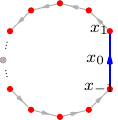}
    \caption{Periodized Loop Graph}
  \end{figure}
\end{minipage}
\begin{minipage}{.65\linewidth}
  \begin{figure}[H]
    \centering
    \includegraphics[width=.8\linewidth]{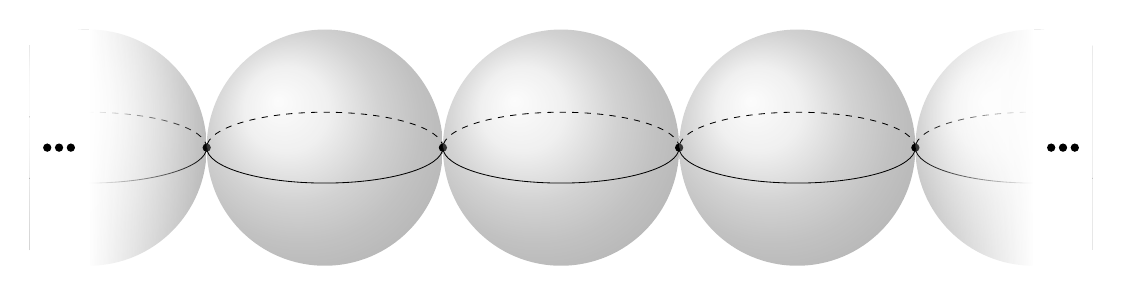}\\
    \includegraphics[]{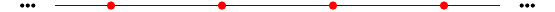}
    \caption{Infinite Chain of $\mathbb{P}^1$}
    \label{fig:inf_chain_P_one}
  \end{figure}  
\end{minipage}

\subsection{Cohomology of Periodized Hypertoric Varieties}

Given $ S \subseteq E $, let
\begin{equation*}
    P_S = \bigcap_{e \in S} H_e^+ \cap \bigcap_{e \notin S} H_e^- \subseteq \mathfrak{t}_\Gamma^*
\end{equation*}
be a polytope bounded by hyperplanes in the hyperplane arrangement $\Hfin$ and
\begin{equation*}
    \toricP[P_S] = \{ [z, w] \in \Mfin \mid z_e = 0 \text{ if } e \in S \text{ and } w_e = 0 \text{ if } e \notin S \}.
\end{equation*}
Then, $\toricP[P_S]$ is a subvariety of $\Mfin$ isomorphic to the toric variety specified by the polytope $P_S$. The core $ \corefin \subseteq \Mfin $ of the hypertoric variety is defined to be the union of $\toricP[P_S]$ over all the bounded polytopes $P_S$. It is proven in \cite[Theorem 6.5]{bielawski2000geometry} that $\Mfin$ admits a $T_\Gamma$-equivariant deformation retraction to $\corefin$. This holds true for $\Mper$. The remaining of this subsection will be devoted to calculating $\HH^\bullet(\coreper)=\HH^\bullet(\Mper)$.
\minormissing{Mper deformation retract to core}

It is well-known that the $0$-dimensional intersections of the hyperplane arrangement $\Hfinchar{\chi}$ corresponds to the fixed points of the $T_\Gamma$-action of $\Mfinchar{\chi, 0}$. When $\chi$ is generic, then the former are labeled by $\face[d]$, where $d$ is the genus of $\Gamma$. Given $ B \subseteq \face[d] $, we let $\fixptP$ be the set of polytopes whose vertices are contained in $B$ and $ \toricB=\bigcup_{P\in \fixptP}\toricP \subseteq \Mfin $ be the union of the toric varieties corresponding to the polytopes in $\fixptP$.

\begin{lemma}[lem:GKM_toric]{}
    Let $ B \subseteq \face[d] $. Suppose that $\fixptP$ is connected and for any $ S \in \face[d-1] $, $ \fixptP \cap \bigcap_{e \in S} H_e $ is either empty or connected. If the odd cohomology of $\toricB$ vanishes, then
    \begin{equation*}
        \HH_{T_\Gamma}^\bullet(\toricB, \mathbb{C}) \cong \frac{\mathbb{C}[E]}{\bigcap_{T^* \in B} \langle e \mid e \in E \setminus T^* \rangle}.
    \end{equation*}
    \noskipline
    \begin{proof}
        Note that the vanishing of odd cohomology implies that $\toricB$ is equivariantly formal. It is also clear that there are only finitely many fixed-points and 1-dimensional orbits. It follows that the GKM conditions are satisfied. Hence, the GKM description \cite[Theorem 1.2.2]{goresky1998equivariant} can be applied to calculate the equivariant cohomology. Note that the $T_\Gamma$-fixed-points of $\toricB$ are indexed by $B$, and two fixed-points are connected by a $T_\Gamma$-invariant sphere if and only if the corresponding points in $\fixptP$ are connected by a line segment in $\fixptP$. If $ S \cup e, S \cup e' \in B $ are of this case, the character of the sphere is given by the unique (up to sign) cycle in $\Gamma \del S$, which equivalently is $ \fundcyclewithtree{e}[\Gamma][E \setminus (S \cup e)] = \pm \fundcyclewithtree{e'}[\Gamma][E \setminus (S \cup e')] $.
        
        With the information above, the GKM description tells us that $\HH_{T_\Gamma}^\bullet(\toricB, \mathbb{C})$ is the subring of $\mathbb{C}[\mathfrak{t}_\Gamma]^B$ consisting of polynomials $(f_{T^*})_{T^* \in B}$ so that for any $ S \cup e, S \cup e' \in B $ connected by a line segment in $\fixptP$, we have
        \begin{equation*}
            f_{S \cup e} - f_{S \cup e'} \in \langle \fundcyclewithtree{e}[\Gamma][E \setminus (S \cup e)] \rangle.
        \end{equation*}
        Since every spanning cotree induces a basis on $\mathfrak{t}_\Gamma$, there is a natural identification 
        \begin{equation*}
            \mathbb{C}[\mathfrak{t}_\Gamma]^B \cong \bigoplus_{T^* \in B} \mathbb{C}[T^*].
        \end{equation*}
        Under this identification, the above condition can be rephrased as $ f_{S \cup e}|_{e=0} = f_{S \cup e'}|_{e'=0} $.

        Now, consider the ring homomorphism $ \mathbb{C}[E] \to \bigoplus_{T^* \in B} \mathbb{C}[T^*] $ which is the quotient map for each summand. The kernel is given by $\bigcap_{T^* \in B} \langle e \mid e \in T \rangle$ and the image is given by
        \begin{equation*}
            \{ (f_{T^*})_{T^* \in B} \mid f_{T^*_1}|_{T^*_1 \setminus T^*_2 = 0} = f_{T^*_2}|_{T^*_2 \setminus T^*_1 = 0} \quad \forall\ T^*_1, T^*_2 \in B \}.
        \end{equation*}
        Here, $f_{T^*}|_{S=0}$ denotes the image of $f_{T^*}$ under the quotient map $ \mathbb{C}[T^*] \to \mathbb{C}[T^* \setminus S] $. Our assumptions on $\fixptP$ guarantee that any $ T^*_1, T^*_2 \in B $ is connected by a path in $\fixptP$ that is contained in $\bigcap_{e \in T^*_1 \cap T^*_2} H_e$. Thus, the image is precisely the subring picked out by the GKM description. The isomorphism in the lemma thus follows from the first isomorphism theorem.
    \end{proof}
\end{lemma}

\begin{lemma}[lem:GKM_odd_cohom_vanish]{}
    Let $ B, B' \subseteq \face[d] $. Suppose that $\fixptP$, $\fixptP[B']$, and $\fixptP[B \cap B']$ all satisfy the conditions in the previous lemma. If the odd cohomology of $\toricB$, $\toricB[B']$ and $\toricB[B \cap B']$ all vanish, then the odd cohomology of $\toricB[B \cup B']$ also vanishes.
    \begin{proof}
        By the previous calculation,
        \begin{equation*}
            \HH^\bullet(\toricB, \mathbb{C}) \cong \frac{\mathbb{C}[E]}{\bigcap_{T^* \in B} \langle e \mid e \in E \setminus T^* \rangle, \langle \HH_1(\Gamma) \rangle},
        \end{equation*}
        and similarly for $\toricB[B']$ and $\toricB[B \cap B']$. Note that the pullback in the cohomology of the inclusion $ \toricB[B \cap B'] \subseteq \toricB $ is given by the natural inclusion of ideals, hence is surjective. The Mayer-Vietoris sequence
        \begin{equation*}
            \begin{tikzcd}
                \cdots \arrow[r] & \HH^\bullet(\toricB[B \cup B'], \mathbb{C}) \arrow[r] & \HH^\bullet(\toricB, \mathbb{C}) \oplus \HH^\bullet(\toricB[B'], \mathbb{C}) \arrow[r] & \HH^\bullet(\toricB[B \cap B'], \mathbb{C}) \arrow[r] & \cdots
            \end{tikzcd}
        \end{equation*}
        then implies that the odd cohomology of $\HH^\bullet(\toricB[B \cup B'], \mathbb{C})$ vanishes.
    \end{proof}
\end{lemma}

For a periodized hyperplane arrangement $\Hper$, the $0$-dimensional intersections of the hyperplane arrangement corresponds to $ \faceper[d] = \{ T^*_I \mid T^* \in \face[d], I \in \ZZ^{T^*} \} $. Given any finite subset $ B \subseteq \faceper[d] $ satisfying the above conditions, we can still apply \Cref{lem:GKM_toric} and \Cref{lem:GKM_odd_cohom_vanish} with $E$ replaced with $E_\per$.
\minormissing{why}

\begin{lemma}[lem:Mper_cohom]{}
    We have
    \begin{equation*}
        \HH_{T_\Gamma}^\bullet(\Mper, \mathbb{C}) \cong \HH_{T_\Gamma}^\bullet(\coreper, \mathbb{C}) \cong \frac{\widehat{\mathbb{C}[E_\per]}}{\bigcap_{T_I^* \in \faceper[d]} \langle e \in E_\per \setminus T_I^* \rangle},
    \end{equation*}
    where $\widehat{\mathbb{C}[E_\per]}$ is a certain completion of $\mathbb{C}[E_\per]$. As an abelian group, this is the infinite product of $\mathbb{C}$ indexed by monomials in $\mathbb{C}[E_\per]$ whose support lies in $\faceper$.
    \begin{proof}
      Pick $T^*=\CB{e_0,\dots ,e_{d-1}}\in \face[d]$ and let $B_{-1}$ be the set of $0$-dimensional intersections of the hyperplane arrangement that lie inside the $d$-dimensional region $P_{-1}$ bounded by $\CB{H_{e_i,0},H_{e_i,1}\mid i\in 0,\dots ,d-1}$.
      Any hyperplane $H$ that passes through $P_{-1}$ will cut it into two sections $B$ and $B'$ with intersections lying inside $H$.
      Inductively, we can glue together first all the polytopes in either section, then glue the two along the intersection. Assuming both $B$ and $B'$ satisfy the assumptions in \Cref{lem:GKM_toric}. Since the intersection $B\cap B'$ contains all $0$-dimensional intersections that lie in $H$, it also satisfies the assumptions in \Cref{lem:GKM_toric}. The base case is when $B$ is precisely the vertices of a polytope. This gives rise to the projective toric variety $\toricB$ with zero odd cohomology. Hence we can apply \Cref{lem:GKM_odd_cohom_vanish} and conclude that the $\toricB[B_{-1}]$ has zero odd cohomology.

      Define $B_{i=(2d)j+(d)k+l}$, where $k=0,1$ and $l<d$, be the set of $0$-dimensional intersections of the hyperplane arrangement that lie inside the region $P_i$ bounded by $\CB{H_{e_j,-d-a_j},H_{e,d+1+b_j}\mid =0,\dots ,d-1}$, where 
      \[
        a_i=\begin{cases}
          1&\text{if }k=1\text{ and }l\geq d,\\
          0&\text{otherwise};
        \end{cases}
        \quad \text{ and } \quad
        b_i=\begin{cases}
          1&\text{if }k=1\text{ or }l\geq d,\\
          0&\text{otherwise}.
        \end{cases}
      \]
      This construction ensures that $P_i$ and $P_{i-1}$ differ by a ``strip" and $P_{i-1}\cap (P_i-P_{i-1})$ lies in a hyperplane. The gluing between $P_{i-1}$ and $P_i-P_{i-1}$ satisfies the assumptions in \Cref{lem:GKM_odd_cohom_vanish}. Hence the odd cohomology of $\toricB[B_i]$ vanishes for all $i$.
      Obviously, we have $ B_{-1} \subseteq B_0 \subseteq \dots \subseteq \faceper[d] $ so that $ \bigcup_{n \in \mathbb{N}} B_n = \faceper[d] $. Equivalently,
        \begin{equation*}
            \varinjlim_{n \in \mathbb{N}} \toricB[B_n] \cong \bigcup_{n \in \mathbb{N}} \toricB[B_n] = \coreper.
        \end{equation*}
        
        By the presentation in \Cref{lem:GKM_toric}, $ \HH_{T_\Gamma}^\bullet(\toricB[B_{n+1}], \mathbb{C}) \to \HH_{T_\Gamma}^\bullet(\toricB[B_n], \mathbb{C}) $ induced by the inclusion $ \toricB[B_n] \subseteq \toricB[B_{n+1}] $ is sujective. Therefore, the inverse system $(\HH_{T_\Gamma}^\bullet(\toricB[B_n], \mathbb{C}))_{n \in \mathbb{N}}$ satisfies the Mittag-Leffler condition. It follows by \cite[Corollary 4.2.4]{kochman1996bordism} that the comparison map
        \begin{equation*}
            \HH_{T_\Gamma}^\bullet(\coreper, \mathbb{C}) \to \varprojlim_{n \in \mathbb{N}} \HH_{T_\Gamma}^\bullet(\toricB[B_n], \mathbb{C}) \cong \varprojlim_{n \in \mathbb{N}} \frac{\mathbb{C}[E_\per]}{\bigcap_{T_I^* \in B_n} \langle e \in E_\per \setminus T^*_I \rangle}
        \end{equation*}
        is an isomorphism. The inverse limit on the right can be described as follows. Let $\widehat{\mathbb{C}[E_\per]}$ be the completion of $\mathbb{C}[E_\per]$ over the direct system of ideals generated by cofinite subsets of $E_\per$. As an abelian group, $\widehat{\mathbb{C}[E_\per]}$ is the infinite product of $\mathbb{C}$ indexed by monomials in $\mathbb{C}[E_\per]$. Since each term in the inverse system is just a quotient of a polynomial ring by a monomial ideal, the inverse limit could be computed, as an abelian group, by taking the infinite product of $\mathbb{C}$ indexed by the union of the monomial bases in the system. Finally, we describe this union. Let $E_\per^\sigma$, where $ \sigma : E_\per \to \mathbb{N} $ is finitely supported, be a monomial in $\mathbb{C}[E_\per]$. Note that it is in the union of the monomial bases, i.e. $ E_\per^\sigma \notin \bigcap_{T_I^* \in \faceper[d]} \langle e \in E_\per \setminus T_I^* \rangle $, if and only if $ \supp(\sigma) \in T_I^* $ for some $ T_I^* \in \faceper[d] $, which is true if and only if $ \supp(\sigma) \in \faceper $. The description in the lemma follows.
    \end{proof}
\end{lemma}

We propose another presentation of the above equivariant cohomology ring. The following construction will illustrate how we ``periodize'' various algebraic constructions later on.
    
\begin{lemma}[lem:SRper]{}
    The short exact sequence in \Cref{def:SRfin} gives a short exact sequence of inverse systems.
    \begin{equation*}
        \begin{tikzcd}
            0 \arrow[r] & (\idealI[\Gamma_n])_{n \in \NN} \arrow[r] & (\ZZ[E_n])_{n \in \NN} \arrow[r] & (\SR[\bullet][\Gamma_n])_{n \in \NN} \arrow[r] & 0
        \end{tikzcd}
    \end{equation*}
    The homomorphisms in all the inverse systems above are surjective. Taking inverse limit in the category of graded abelian groups gives a short exact sequence as follows.
    \begin{equation*}
        \begin{tikzcd}
            0 \arrow[r] & \idealIper \arrow[r] & \varprojlim_{n \in \NN} \ZZ[E_n] \arrow[r] & \SRper \arrow[r] & 0
        \end{tikzcd}
    \end{equation*}
    In particular, $\SRper[2k]$ is the infinite product of $\ZZ$ indexed by degree $k$ monomials in $\ZZ[E_\per]$ whose support lies in $\faceper$. Moreover, the completed tensor product $\SRper \mathbin{\widehat{\otimes}}_\ZZ \mathbb{C}$ coincides with the description of $\HH_{T_\Gamma}^\bullet(\Mper, \mathbb{C})$ in \Cref{lem:Mper_cohom}.
    \begin{proof}
        Given $ \sigma : E_n \to \NN $, denote by $E_n^\sigma$ the monomial $ \prod_{e_i \in E_n} e_i^{\sigma(e_i)} $ in $\ZZ[E_n]$. The contraction map $ \ZZ[E_{n+1}] \to \ZZ[E_n] $ sends $(E_{n+1})^\sigma$ to $E_n^{\sigma|_{E_n}}$ if $ \supp(\sigma) \subseteq E_n $, and $0$ else. It is clear that these maps commute with the maps in a short exact sequence. We compute $\SRper$ directly. The other inverse limits are computed analogously. Note that $\SR[2k][\Gamma_n]$ as an abelian group is given by the set of all functions
        \begin{equation*}
            \left\{ E_n^\sigma \mid \sum_{e_i \in E_n} \sigma(e_i) = k \text{ and } \supp(\sigma) \in \face[\bullet][\Gamma_n] \right\} \to \ZZ.
        \end{equation*}
        Thus, its inverse limit is given by the set of all functions
        \begin{equation*}
            \varinjlim_{n \in \NN} \left\{ E_n^\sigma \mid \sum_{e_i \in E_n} \sigma(e_i) = k \text{ and } \supp(\sigma) \in \face[\bullet][\Gamma_n] \right\} \cong \left\{ E_\per^\sigma \mid \sum_{e_i \in E_\per} \sigma(e_i) = k \text{ and } \supp(\sigma) \in \faceper \right\} \to \ZZ.
        \end{equation*}
        This clearly coincides with the description in \Cref{lem:Mper_cohom}.
    \end{proof}
\end{lemma}

We suspect that the integral equivariant cohomology of $\Mper$ is precisely $\SRper$, though we do not have a proof for that. Still, we will work with $\SRper$ over $\ZZ$ instead of over $\mathbb{C}$. Moreover, the algebraic arguments in the rest of the paper work over both $\ZZ$ and $\mathbb{C}$, and does not depend on the actual cohomology of $\Mper$.
\minormissing{integral cohom of Rper}

\begin{lemma}[lem:Rper]{}
    There is an exact sequence as follows periodizing the one in \Cref{def:Rfin}.
    \begin{equation*}
        \begin{tikzcd}
            \SRper \otimes \HH_1(\Gamma) \arrow[r] & \SRper \arrow[r] & \Rper \arrow[r] & 0
        \end{tikzcd}
    \end{equation*}
    Here, the map $ \HH_1(\Gamma) \to \SRper $ is the restriction of the map $ \ZZ^E \to \SRper $ given by $ e \mapsto \sum_{i \in \ZZ} e_i $. Moreover, the completed tensor product $\Rper \mathbin{\widehat{\otimes}}_\ZZ \mathbb{C}$ coincides with $\HH^\bullet(\Mper, \mathbb{C})$.
    \begin{proof}
        Note that $ \HH_1(\Gamma_n) \cong \HH_1(\Gamma) $ for all $ n \in \NN $. Moreover, the periodization of the map $ \HH_1(\Gamma) \hookrightarrow \ZZ^E $ is the map $ \HH_1(\Gamma) \to \ZZ^{E_\per} $ given by $ \gamma \mapsto \sum_{e \in E} \pair{\gamma}{e} \left(\sum_{i \in \ZZ} e_i\right) $. The description of $ \HH_1(\Gamma) \to \SRper $ thus follows.
    \end{proof}
\end{lemma}

\subsection{Periodization for \texorpdfstring{$\HTfin$}{HTfin}}
\label{sec:periodization_HTfin}

We will now periodize $\HTfin$ and its homotopy equivalence with $\Rfin$. By our description of $\face[\bullet][\Gamma_n]$, it is clear that
\begin{equation*}
    \HTfin[\bullet][\bullet][\Gamma_n] = \bigoplus_{S_I \in \face[\bullet][\Gamma_n]} \Ext \HH_1(\Gamma_n \del S_I) \cong \bigoplus_{S \in \face} \ZZ^{[-n, n]^S} \otimes \Ext \HH_1(\Gamma \del S).
\end{equation*}
The contraction map $ \HTfin[\bullet][\bullet][\Gamma_{n+1}] \to \HTfin[\bullet][\bullet][\Gamma_n] $ is induced by the inclusion $ [-n, n] \hookrightarrow [-n-1, n+1] $ of sets. To see that the contraction maps commutes with the differentials, it suffices to check that for $ S_I \in \face[\bullet][\Gamma_{n+1}] $ with $ \im I \subseteq [-n, n] $, the following diagram commutes.
\begin{equation*}
    \begin{tikzcd}
        \HH_1(\Gamma_{n+1} \del S_I) \arrow[r, hook] \arrow[d] & \ZZ^{E_{n+1} \setminus S_I} \arrow[d] \\
        \HH_1(\Gamma_n \del S_I) \arrow[r, hook] & \ZZ^{E_n \setminus S_I}
    \end{tikzcd}
\end{equation*}
This is true by the formula of the bottom arrow $ \gamma \mapsto \sum_{e \in E \setminus S} \pair{\gamma}{e}[\Gamma \del S] \left(\sum_{i \in [-n, n]} e_i\right) $. It follows that we have a periodization of $\HTfin$, which we explicitly define as follows.

\begin{definition}[def:]{}
    The complex $(\HTfin, \dHTfin)$ periodizes to $(\HT, \dHT)$ given by
    \begin{equation*}
        \HT[2p][q] = \prod_{S_I \in \faceper[p]} \Ext[q] \HH_1(\Gamma \del S) = \bigoplus_{S \in \face[p]} \ZZ^{\ZZ^S} \otimes \Ext[q] \HH_1(\Gamma \del S),
    \end{equation*}
    where the differential $ \dHT : \HT \to \HT[\bullet+2][\bullet-1] $ restricted to the summand $\ZZ^{\ZZ^S} \otimes \Ext \HH_1(\Gamma \del S) \to \ZZ^{\ZZ^{S \cup e}} \otimes \Ext[\bullet-1] \HH_1(\Gamma \del S \del e)$ is induced by
    \begin{itemize}
        \item the coordinate projection $ \ZZ^{S \cup e} \to \ZZ^S $, and
        \item the interior product $ \intp{e} : \Ext \HH_1(\Gamma \del S) \to \Ext[\bullet-1] \HH_1(\Gamma \del S \del e) $ with $[e] \in \HH^1(\Gamma \del S)$.
    \end{itemize}

    Define $ \grHT[k][p] = \HT[2p][k-p] $. Then, each $\grHT[k]$ is a subcomplex with differential of degree $1$, giving a decomposition $ \HT \cong \bigoplus_{k \in \mathbb{N}} \grHT[k]{[-k]} $.
\end{definition}

Analogous to the finite case, we would expect the following homotopy equivalence.
\begin{theorem}[thm:HT_eqv]{}
    There are homotopy equivalences
    \begin{equation*}
        \begin{tikzcd}
            (\grHT[k], \dHT) \arrow[r, "\fHT", shift left] & \arrow[l, "\gHT", shift left] \Rper[2k]{[-k]}.
        \end{tikzcd}
    \end{equation*}
    \noskipline
    \begin{remark}
        In particular, the cohomology of $\HT$, graded by the total degree, is isomorphic to $\Rper$ as graded abelian groups.
    \end{remark}
\end{theorem}

The construction of the above homotopy equivalence depends on the original choices in \Cref{thm:HTfin_eqv} and a further choice as indicated below. To show this, it suffices to show that $\fHTfin$, $\gHTfin$ and $\hHTfin$ commute with the contraction maps, and then pass them to the inverse limit. This is clear for $\fHTfin$. For $\gHTfin$, recall that a coherent cotree on $\Gamma$ defines a natural coherent cotree on $\Gamma_n$ for all $ n \in \NN $. Since $\gHTfin$ is uniquely determined by the split in \Cref{thm:Rfin_split}, it suffices to show that the split is compatible with the contractions, i.e. $\ZZ^{\basis[\bullet][\Gamma_{n+1}]}$ is mapped to $\ZZ^{\basis[\bullet][\Gamma_n]}$ under contraction. Recall our description of $\basis[\bullet][\Gamma_n]$ in \Cref{lem:cotreein_periodization}. Let $ S_I \in \basis[\bullet][\Gamma_{n+1}] $. Then, $ \cotreein{S} \subseteq \supp(I) $. It follows that $ S_I \in \basis[\bullet][\Gamma_n] $ if $ \im I \subseteq [-n, n]^S $. To define $\hHTfin$ for $\Gamma_n$, we first need to fix a choice $ \choice{S_I}_{\Gamma_n} \in \cotreein{S_I}[\Gamma_n] $ for all $ S_I \in \face[\bullet][\Gamma_n] \setminus \basis[\bullet][\Gamma_n] $. For all $ S \in \face $ and $ S' \subseteq S $ such that $ S \setminus S' \in \face[\bullet][\Gamma \del S'] \setminus \basis[\bullet][\Gamma \del S'] $, pick $ \choice{S, S'} \in \cotreein{S \setminus S'}[\Gamma \del S'] = \cotreein{S} \cap (S \setminus S') $. Then for $ S_I \in \face[\bullet][\Gamma_n] \setminus \basis[\bullet][\Gamma_n] $, define $ \choice{S_I}_{\Gamma_n} = e_0 $, where $ e = \choice{S, \supp(I)} $.
\begin{lemma}[lem:]{}
    The homotopies $\hHTfin_{\Gamma_n}$ defined by the choices $\choice{-}_{\Gamma_n}$ commute with the contraction maps.
    \begin{proof}
        The homotopy $\hHTfin_{\Gamma_n}$ maps $1_{S_I} x^1 \dots x^k$ to $0$ if either $k=0$ and $ S \setminus \supp(I) \in \basis[\bullet][\Gamma \del \supp(I)] $, or $\choice{S \cup x^1 \dots x^k, \supp(I)} \notin S $. Otherwise, if $ x^0 = \choice{S \cup x^1 \dots x^k, \supp(I)} \in S $, then it is mapped to
        \begin{equation*}
            1_{S_I \setminus x^0_0} x^0 \dots x^k - \sum_{t_i \in \tree{S_I \setminus x^0_0}[\Gamma_n]} \hHTfin_{\Gamma_n}\left(1_{S_I \setminus x^0_0 \cup t_i} \intp{t}(x^0 \dots x^k)\right).
        \end{equation*}
        Let $ 1_{S_I} x^1 \dots x^k \in \HTfin[\bullet][\Gamma_{n+1}] $. Suppose $ \im I \cap \{ -n-1, n+1 \} \neq \emptyset $. Note that this property is preserved by $\hHTfin_{\Gamma_{n+1}}$. It follows that this element, both before and after applying $\hHTfin_{\Gamma_{n+1}}$, is mapped to $0$ under the contraction map. Now suppose $ \im I \subseteq [-n, n] $. In this case, the contraction map is the identity and the description of $\hHTfin_{\Gamma_n}$ and $\hHTfin_{\Gamma_{n+1}}$ are identical. It follows that the homotopies commute with the contraction maps.
    \end{proof}
\end{lemma}

Analogous to \Cref{cor:Rfin_les}, we have the following corollary.
\begin{corollary}[cor:Rfper_les]{}
    The homotopy equivalence gives a long exact sequence as follows.
    \begin{equation*}
        \begin{tikzcd}
            0 \arrow[r] & \grHT[k][0] \arrow[r, shift left, "\dHT"] & \cdots \arrow[l, shift left, dashed, "\hHT"] \arrow[r, shift left, "\dHT"] & \grHT[k][k-1] \arrow[l, shift left, dashed, "\hHT"] \arrow[r, shift left, "\dHT"] & \grHT[k][k] \arrow[l, shift left, dashed, "\hHT"] \arrow[r, shift left, "\fHT"] & \Rper[2k] \arrow[r] \arrow[l, shift left, dashed, "\gHT"] & 0
        \end{tikzcd}
    \end{equation*}
\end{corollary}

\subsection{Deletion-Contraction for Periodized Hypertoric Varieties}

We now periodized the deletion-contraction short exact sequence in \Cref{thm:Rfin_delcon} by starting with the following short exact sequence of polynomial rings.
\begin{equation*}
    \begin{tikzcd}
        0 \arrow[r] & \ZZ[(E \setminus e)_n]^{[-n, n]} \arrow[r] & \ZZ[E_n] \arrow[r] & \ZZ[(E \con e)_n] \arrow[r] & 0
    \end{tikzcd}
\end{equation*}
The inclusion is given by $ (r_i)_{i \in [-n, n]} \mapsto \sum_{i \in [-n, n]} r_ie_i $ and the projection is the quotient map realizing the isomorphism $ \ZZ[(E \setminus e)_n] \cong \ZZ[E_n] / \langle e_i \mid i \in [-n, n] \rangle $.

\begin{theorem}[thm:]{}
    If $ e \in E $ is neither a bond nor a loop, there is a deletion-contraction short exact sequence
    \begin{equation*}
        \begin{tikzcd}
            0 \arrow[r] & \Rfin[\bullet-2][(\Gamma \del e)_n]^{[-n, n]} \arrow[r] & \Rfin[\bullet][\Gamma_n] \arrow[r] & \Rfin[\bullet][(\Gamma \con e)_n] \arrow[r] & 0.
        \end{tikzcd}
    \end{equation*}
    \noskipline
    \begin{proof}
        The proof that the maps descend to the quotients is similar. We will only show that this sequence is short exact. Again, fix a coherent cotree $\cotreefunction$ such that $ e \in \tree{\emptyset} $, which gives $\cotreefunction[(\Gamma \del e)_n]$, $\cotreefunction[\Gamma_n]$ and $\cotreefunction[(\Gamma \con e)_n]$ as constructed in \Cref{sec:periodization_of_graphs}. By \Cref{lem:cotreein_periodization},
        \begin{equation*}
            \basis[\bullet][\Gamma_n] = \{ S_I \subseteq E_n \mid \cotreein{S} \subseteq \supp(I) \}.
        \end{equation*}
        For deletion, we have
        \begin{equation*}
            \basis[\bullet][(\Gamma \del e)_n] = \{ S_I \subseteq (E \setminus e)_n \mid \cotreein{S}[\Gamma \del e] \subseteq \supp(I) \}.
        \end{equation*}
        By \Cref{lem:cotreein_subset}, $ \cotreein{S}[\Gamma \del e] = \cotreein{S \cup e} $. Thus, for each $ i \in [-n, n] $, there is an inclusion $ \basis[\bullet][(\Gamma \del e)_n] \hookrightarrow \basis[\bullet][\Gamma_n] $ given by $ S_I \mapsto S_I \cup e_i $. For contraction, we have
        \begin{equation*}
            \basis[\bullet][(\Gamma \con e)_n] = \{ S_I \subseteq (E \setminus e)_n \mid \cotreein{S}[\Gamma \con e] \subseteq \supp(I) \}.
        \end{equation*}
        It is easy to see that $ \cotreein{S}[\Gamma \con e] = \cotreein{S} $. Thus, we also have an inclusion $ \basis[\bullet][(\Gamma \con e)_n] \subseteq \basis[\bullet][\Gamma_n] $. These inclusions together give a partition of $\basis[\bullet][\Gamma_n]$, proving the exactness.
    \end{proof}
\end{theorem}

It is easy to see that this short exact sequence commutes with the contraction maps $ \Gamma_{n+1} \to \Gamma_n $. Since the inverse systems of rings $\{\Rfin[\bullet][\Gamma_n]\}_{n \in \mathbb{N}}$ satisfies the Mittag-Leffler condition. We thus obtain the following by taking inverse limit in the category of graded abelian groups.
\begin{theorem}[thm:Rper_delcon]{}
    If $ e \in E $ is neither a bond nor a loop, there is a deletion-contraction short exact sequence as follows.
    \begin{equation*}
        \begin{tikzcd}
            0 \arrow[r] & \Rper[\bullet-2][\Gamma\del e]^{\mathbb{Z}} \arrow[r] & \Rper \arrow[r] & \Rper[\bullet][\Gamma\con e] \arrow[r] & 0
        \end{tikzcd}
      \end{equation*}
\end{theorem}

%% file: sec_multiplicative_hypertoric_and_cks.tex
\section{Multiplicative Hypertoric Variety and the CKS Complex}
\label{sec:multiplicative_hypertoric_and_cks}

\subsection{Lattice Action on Periodized Hypertoric Varieties}

We now describe a lattice action on the periodized hypertoric varieties.

Consider the $\ZZ$-action on $\Mper[\loopgraph]$ given by translations $ [i, z, w] \mapsto [i+1, z, w] $. This action is compatible with the $\mathbb{C}^\times$-action and the moment map $\mu^\per$.

This extends to a $\mathfrak{d}_{\Gamma, \ZZ}^*$-action on $\Mper[\loopgraph]^E$. This action preserves the moment map condition, while the stability condition is only preserved by the action of the sublattice $ \mathfrak{t}_{\Gamma, \ZZ}^* \subseteq \mathfrak{d}_{\Gamma, \ZZ}^* $, since it simultaneously translates all the hyperplanes. Thus, there is a $\mathfrak{t}_{\Gamma, \ZZ}^* = \HH_1(\Gamma)$-action on the quotient $\Mperchar{\chi, \lambda}$.

There is a corresponding $\HH_1(\Gamma)$-action on $\Rper$. For any $ S \in \face $, there is an action of $\ZZ^S$ on $\ZZ^{\ZZ^S}$ induced by the translation action on $\ZZ^S$. Explicitly, the action is given by
\begin{equation*}
    \delta \cdot (a_I)_{I \in \ZZ^S} = (a_{I+\delta})_{I \in \ZZ^S}
\end{equation*}
for all $ \delta \in \ZZ^S $ and $ a = (a_I)_{I \in \ZZ^S} \in \ZZ^{\ZZ^S} $. We extend this to an action by $\ZZ^E$ along the projection $ \ZZ^E \to \ZZ^S $. This gives an action of $\ZZ^E$ on $\HT$, which commutes with the differential since the translation action of $\ZZ^E$ on the abelian groups $\{ \ZZ^S \}_{S \in \face}$ commutes with the coordinate projection $ \ZZ^{S'} \to \ZZ^S $ for any $ S \subseteq S' $. Therefore, this descends to a $\ZZ^E$-action on $\Rper$ along $\fHT$, which is further restricted to a $\HH_1(\Gamma)$-action. From the GKM perspective, this action corresponds to the $\HH_1(\Gamma)$-action on the hyperplanes in $\Hperchar{\chi}$ by $ \gamma \cdot H_{e, i} \mapsto H_{e, i + \pair{\gamma}{e}} $, and thereby the $0$-dimensional intersections of the hyperplanes and the corresponding torus fixed-points of $\Mperchar{\chi, 0}$.
\minormissing{The lattice action on $\Mper$ corresponds to that in $\Rper$}

Identifying $\Rper[\bullet-2][\Gamma\del e]^{\mathbb{Z}}$ as an ideal in $\Rper$, and noting that $ \HH_1(\Gamma \con e) \cong \HH_1(\Gamma) $, it is not hard to see that the short exact sequence in \Cref{thm:Rper_delcon} is $\HH_1(\Gamma)$-equivariant.

\subsection{Multiplicative Hypertoric Variety}

The $\ZZ$-action on $\Mper[\loopgraph]$ is not free. For example, the point $[0, 1, e^{2\pi i/n}]$ has stablizer $n\ZZ$. However, we have the following result.

\begin{proposition}[prop:]{\cite[Proposition 7.1.2]{dancso2024deletion}}
    The $\ZZ$-action on $\Mper[\loopgraph]$ is free and properly discontinuous on $(\mu^\per)^{-1}(\opendisc)$, where $ \opendisc = \{ z \in \mathbb{C} \mid |z| < 1 \} $.
\end{proposition}

In fact, \cite{dancso2024deletion} showed the existance of a map $ (\mu^\per)^{-1}(\opendisc) \to \mathbb{R} $ that is both $\ZZ$-equivariant and $\mathbb{U}_1$-equivariant, where $\mathbb{U}_1$ denotes the maximal compact torus in $\mathbb{C}^\times$. This extends to a map from $ (\nu^\per_\Gamma)^{-1}(\mathbb{D}^E \cap \mathfrak{t}_{\Gamma}^*) $ to $\mathfrak{t}_{\Gamma, \mathbb{R}}^*$ that is both $\HH_1(\Gamma)$-equivariant and $T_{\Gamma, \mathrm{cpt}}$-equivariant, where $ \nu^\per_\Gamma : \Mperchar{\chi, 0} \to \mathfrak{t}_{\Gamma}^* $ is the complex moment map for the $T_\Gamma$-action on $\Mperchar{\chi, \lambda}$, and $\mathbb{D}^E$ is the unit polydisc in $\mathfrak{d}_\Gamma^*$. It follows that the $\HH_1(\Gamma)$-action on $(\nu^\per_\Gamma)^{-1}(\mathbb{D}^E \cap \mathfrak{t}_{\Gamma}^*)$ again free and properly discontinuous.

\minormissing{lattice action free and discontinuous on preimage of disc}

We define the multiplicative hypertoric variety $\Mmul$ to be the quotient of $ (\nu^\per_\Gamma)^{-1}(\mathbb{D}^E \cap \mathfrak{t}_{\Gamma}^*) \subseteq \Mperchar{\chi, \lambda} $ by the $\HH_1(\Gamma)$-action.

\answer{When defining the multiplicative hypertoric, you need to introduce the notation you're going to use for it (presumably $\Mmul$}


\begin{example}[exp:]{}
    The $\ZZ$-invariant map $ \mu^\per : (\mu^\per)^{-1}(\opendisc) \to \opendisc $ descends to an elliptic fibration $ \Mmul[\loopgraph] \to \opendisc $ with central fiber a nodal elliptic curve. The node corresponds to the $\ZZ$-orbit $ \{ [i, 0, 0] \mid i \in \ZZ \} \subseteq \Mper[\loopgraph] $.
\end{example}

\subsection{The CKS Complex and Deletion-Contraction}
\label{sec:cks_del_con}

It is proven in \cite{dancso2024deletion} that the cohomology of multiplicative hypertoric varieties is isomorphic to the cohomology of the CKS complex.
This complex first appeared in \cite{cattani19872} and later in \cite{migliorini2021support} that \cite{dancso2024deletion} is based on.
We now give the definition of this complex.

\begin{definition}[def:CKS]{}
    The CKS complex $\CKS$ is defined by
    \begin{equation*}
        \CKS[2p][q][r] = \bigoplus_{S \in \face[p]} \Ext[q] \HH_1(\Gamma \del S) \otimes \Ext[r] \HH^1(\Gamma \del S).
    \end{equation*}
    Given $ S \in \face[p] $, $ \gamma_1, \dots, \gamma_q \in \HH_1(\Gamma \del S) $ and $ \alpha_1, \dots, \alpha_r \in \HH^1(\Gamma \del S) $, the corresponding element in $\CKS[2p][q][r]$ is denoted by $ 1_S \gamma_1 \dots \gamma_q \otimes \alpha_1 \dots \alpha_r $. The differential $ \dCKS : \CKS \to \CKS[\bullet+2][\bullet-1][\bullet] $ restricted to the summands
    \begin{equation*}
        \Ext \HH_1(\Gamma \del S) \otimes \Ext \HH^1(\Gamma \del S) \to \Ext[\bullet-1] \HH_1(\Gamma \del S \del e) \otimes \Ext \HH^1(\Gamma \del S \del e)
    \end{equation*}
    is induced from
    \begin{itemize}
        \item the natural projection $ \HH^1(\Gamma \del S) \to \HH^1(\Gamma \del S \del e) $, and
        \item the interior product $ \intp{e} : \Ext \HH_1(\Gamma \del S) \to \Ext[\bullet-1] \HH_1(\Gamma \del S \del e) $ with $ [e] \in \HH^1(\Gamma \del S) $.
    \end{itemize}

    Define $ \grCKS[k][\ell][p] = \CKS[2p][k-p][\ell] $. Then, each $\grCKS[k][\ell]$ is a subcomplex with differential of degree $1$, giving a decomposition $ \CKS \cong \bigoplus_{k, \ell \in \mathbb{N}} \grCKS[k][\ell]{[-k-\ell]} $.
    \begin{remark}
        In \cite{dancso2024deletion}, the CKS complex is bigraded by $(2p, q+r)$. And the cohomological grading is given by $2p+q+r$.
    \end{remark}
\end{definition}

\begin{example}[exp:cks_loopgraph]{}
    For $\Gamma = \loopgraph$, the differential is an isomorphism $ \Ext[1] \HH_1(\Gamma) \to \Ext[0] \HH_1(\Gamma \del e) $ and $0$ elsewhere, where $e$ is the unique edge. It follows that the cohomology is given by $ \ZZ \oplus \HH^1(\Gamma) \oplus (\HH_1(\Gamma) \otimes \HH^1(\Gamma)) $.
    This agrees with the geometric picture that the core of $\Mmul[\loopgraph]$ is a pinched torus.
\end{example}

Let $ e \in E $ be an edge that is neither a bond nor a loop. Using the deletion-contraction sequence for the coindependence complex in \Cref{eqn:face_delcon}, and the fact that $ \HH_1(\Gamma \con e) \cong \HH_1(\Gamma) $ and $ \HH^1(\Gamma \con e) \cong \HH^1(\Gamma) $, it is easy to construct a deletion-contraction sequence for the CKS complex.

\begin{theorem}[thm:cks_del_con]{\cite{dancso2024deletion}}
    If $e\in E$ is neither a bond nor a loop, there is a deletion-contraction short exact sequence of the CKS complex
    \begin{equation*}
        \begin{tikzcd}
            0 \arrow[r] & \CKS[\bullet-2][\bullet][\bullet][\Gamma \del e] \arrow[r] & \CKS \arrow[r] & \CKS[\bullet][\bullet][\bullet][\Gamma \con e] \arrow[r] & 0
        \end{tikzcd}
    \end{equation*}
    which induces a long exact sequence in cohomology.
\end{theorem}

In fact, we can decompose this sequence into finer pieces according to the two extra gradings. This gives short exact sequences
\begin{equation*}
    \begin{tikzcd}
        0 \arrow[r] & \grCKS[k-1][\ell][\bullet-1][\Gamma \del e] \arrow[r] & \grCKS[k][\ell] \arrow[r] & \grCKS[k][\ell][\bullet][\Gamma \con e] \arrow[r] & 0
    \end{tikzcd}
\end{equation*}
of cochain complexes, and also the corresponding long exact sequences in cohomology.

Define
\begin{equation*}
    \euler{k}{\ell} = \chi(\grCKS[k][\ell]) = \sum_{p \in \mathbb{N}} (-1)^p \dim_\ZZ \HH^{2p, k-p, \ell}(\CKS, \dCKS).
\end{equation*}
By a standard homological algebra argument, the deletion-contraction sequence of the CKS complex implies that $ \euler{k}{\ell} = \euler{k}{\ell}[\Gamma \con e] - \euler{k-1}{\ell}[\Gamma \del e] $. Define a polynomial 
\begin{equation*}
    \HTutte{x, y} = (-1)^d \sum_{k, \ell = 0}^d \euler{k}{\ell} \, x^{d-k} y^\ell,
\end{equation*}
where $d$ is the genus of $\Gamma$. In fact, this polynomial can be expressed in terms of the Tutte polynomial of $\Gamma$.

\begin{theorem}[thm:]{}
    $\HTutte{x,y} = \Tutte{-(x+y+xy),1}$.
    \begin{proof}
        By the universality of the Tutte polynomial, any graph invariant satisfying
        \[
            \HTutte{x,y}=
            \begin{cases}
                \HTutte[\Gamma\del e]{x,y}+\HTutte[\Gamma\con e]{x,y} & \textnormal{if $e \in E$ is neither a bond nor a loop}, \text{ and} \\
                \HTutte[\bridgegraph]{x,y}^b\HTutte[\loopgraph]{x,y}^l & \textnormal{if $\Gamma$ consists of exactly $b$ bridges (one-edge bonds) and $l$ loops}.
            \end{cases}
        \]
        is an evaluation of the Tutte polynomial, given by $ \HTutte{x,y} = \Tutte{\HTutte[\loopgraph]{x,y},\HTutte[\bridgegraph]{x,y}} $. The discussion preceding this theorem showed that $\HTutte{x,y}$ satisfies the deletion-contraction condition.

        It is clear that a wedge product of graphs corresponds to the tensor product of complexes. Thus, it suffices to study the above polynomial for the loop graph and the bridge graph. For the former, \Cref{exp:cks_loopgraph} tells us that the cohomology has three generators of degree $(0, 0, 0)$, $(0, 0, 1)$ and $(0, 1, 1)$ respectively, giving $ \euler{0}{0}[\loopgraph] = \euler{0}{1}[\loopgraph] = \euler{1}{1}[\loopgraph] = 1 $ and $ \HTutte[\loopgraph]{x, y} = -(x+y+xy) $. For the latter, it is clear that the cohomology has only one generator of degree $(0, 0, 0)$, giving $ \euler{0}{0}[\bridgegraph] = 1 $ and $ \HTutte[\bridgegraph]{x, y} = 1 $.
    \end{proof}
\end{theorem}


%% file: sec_gcnew.tex
\section{Computing the Group Cohomologies}
\label{sec:gcnew}

\subsection{Group Cohomology}

For any abelian group $G$, denote by $ \invar{G} = \Hom_{\ZZ[G]}(\ZZ, -) $ the invariant submodule functor, where $\ZZ$ has trivial $G$-action. For a $\ZZ[G]$-module $M$, the group cohomology of $G$ with coefficients in $M$ is defined by
\begin{equation*}
    \HH^\bullet_{\mathrm{Grp}}(G, M) = \Rder\invar{G}(M),
\end{equation*}
where $\Rder\invar{G}$ denote the right derived functor of $\invar{G}$. This can be computed as $\HH^\bullet(\Hom_{\ZZ[G]}(P_\bullet, M))$ for any projective resolution $P_\bullet$ of $\ZZ$ as a $\ZZ[G]$-module. We will only be dealing with finitely generated free abelian groups. For the sake of completeness, let us record a standard projective resolution in this case.

Pick a basis $\{ e_1, \dots, e_n \}$ of $G$, where $ n = \dim_\ZZ G $. For $ g \in G $, write the corresponding generator in $\ZZ[G]$ as $1_g$. A projective resolution is given by
\begin{equation*}
    \begin{tikzcd}
        0 \arrow[r] & \Ext[n]_{\ZZ} G \otimes_{\ZZ} \ZZ[G] \arrow[r, "\partial"] & \dots \arrow[r, "\partial"] & \Ext[1]_{\ZZ} G \otimes_{\ZZ} \ZZ[G] \arrow[r, "\partial"] & \Ext[0]_{\ZZ} G \otimes_{\ZZ} \ZZ[G] \arrow[r, "\epsilon"] & \ZZ \arrow[r] & 0,
    \end{tikzcd}
\end{equation*}
where
\begin{equation*}
    \partial : \ZZ[G]^n \cong \Ext[1]_{\ZZ} G \otimes_{\ZZ} \ZZ[G] \to \Ext[0]_{\ZZ} G \otimes_{\ZZ} \ZZ[G] \cong \ZZ[G]
\end{equation*}
is given by $ (1_{e_i} - 1_0)_{i=1}^n \in \ZZ[G]^n \cong \Hom_{\ZZ[G]}(\ZZ[G]^n, \ZZ[G]) $, which is then extended to higher exterior powers via the Leibniz rule, and $\epsilon$ is given by $ 1_g \mapsto 1 $ for all $ g \in G $.

\subsection{Some Homological Algebra Results}

\subsubsection{Degeneration of the Spectral Sequence of Double Complex}

We record a degeneration result about the spectral sequence of double complex. The proof of the following theorem can be found in \cite[Theorem 2.3.4]{neukirch2013cohomology}.

\begin{theorem}[thm:degen_spectral_seq]{}
    Let $(C^\bullet, d_C)$ and $(D^\bullet, d_D)$ be bounded below complexes of abelian groups so that $C^\bullet$ is degreewise flat. Then, the spectral sequence of the double complex $C^\bullet \otimes D^\bullet$ with the total differential $d_{\mathrm{tot}}$ given by
    \begin{equation*}
        E_2^{p, q} = \HH^p(\HH^{\bullet, q}(C^\bullet \otimes D^\bullet, d_D), d_C) \quad \Rightarrow \quad \HH^\bullet(C^\bullet \otimes D^\bullet, d_{\mathrm{tot}})
    \end{equation*}
    degenerates at $E_2$ page. Moreover, there are non-canonical isomorphisms
    \begin{equation*}
        \HH^n(C^\bullet \otimes D^\bullet, d_{\mathrm{tot}}) \cong \bigoplus_{p+q=n} E_2^{p, q}.
    \end{equation*}
\end{theorem}

\subsubsection{Computation of Group Cohomology}

Let
\begin{equation*}
    \begin{tikzcd}
        1 \arrow[r] & N \arrow[r, "\iota"] & G \arrow[r, "\pi"] & H \arrow[r] & 1
    \end{tikzcd}
\end{equation*}
be a short exact sequence of abelian groups. Given any $\ZZ[G]$-modules $A$ and $B$, we define the $\ZZ[G]$-module $B^A$ to be the set of all functions from $A$ to $B$ with a $G$-action by conjugation, i.e.
\begin{equation*}
    (g \cdot f)(a) = g \cdot f(g^{-1} \cdot a)
\end{equation*}
for all $ g \in G $, $ a \in A $ and $ f \in B^A $. Equivalently, $ B^A \cong \Hom_\ZZ(\ZZ^{\oplus A}, B) $ as $\ZZ[G]$-modules. We will show the following result.
\begin{theorem}[thm:isom_gc]{}
    For any $\ZZ[G]$-module $A$, $ \Rder\invar{G}(A^H) \cong \Rder\invar{N}(A) $ as abelian groups.
\end{theorem}

We start with the following easy facts.

\begin{lemma}[lem:isom_G_mod]{}
    Let $A$ be a $\ZZ[G]$-module. Then, $ A^G \cong (A_{\mathrm{triv}})^G $ as $\ZZ[G]$-modules, where $A_{\mathrm{triv}}$ is the underlying abelian group of $A$ with trivial $G$-action.
    \begin{proof}
        Define $ \phi : A^G \to (A_{\mathrm{triv}})^G $ by $ \phi(f)(g) = g^{-1}f(g) $. It is easy to see that this is compatible with the $G$-action on both sides, and is also bijective.
    \end{proof}
\end{lemma}

\begin{lemma}[lem:gc_con_at_zero]{}
    For any $\ZZ[G]$-module $A$, $ \Rder\invar{G}(A^G) \cong A $ as abelian groups.
    \begin{proof}
        By \Cref{lem:isom_G_mod}, we can assume that $A$ has a trivial $G$-action. Then, 
        \[
            \Rder\invar{G}(A^G) = \Rder\Hom_{\ZZ[G]}(\ZZ, \Hom_\ZZ(\ZZ[G], A)).
        \]
        Since $\ZZ[G]$ is a free abelian group, we can replace the second hom functor with a derived one. By the derived tensor-hom adjunction,
        \begin{equation*}
            \mathbb{R}\Hom_{\ZZ[G]}(\ZZ, \mathbb{R}\Hom_\ZZ(\ZZ[G], A)) \cong \mathbb{R}\Hom_\ZZ(\ZZ[G] \otimes_{\ZZ[G]}^{\mathbb{L}} \ZZ, A) \cong \mathbb{R}\Hom_\ZZ(\ZZ, A) \cong A.
        \end{equation*}
    \end{proof}
\end{lemma}

Given any $\ZZ[G]$-module $A$, we have
\begin{equation*}
    \Hom_{\ZZ[G]}(\ZZ, A) \cong \Hom_{\ZZ[H]}(\ZZ, \Hom_{\ZZ[N]}(\ZZ, \mathrm{Res}_{\ZZ[G]}^{\ZZ[N]}(A)))
\end{equation*}
as a consequence of the tensor-hom adjunction, where all $\ZZ$ above have trivial group actions. By abuse of notation, we record the above as $ \invar{G} \cong \invar{H} \circ \invar{N} $. Therefore, we have the Grothendieck spectral sequence
\begin{equation*}
    E_2^{p, q} = \Rder[p]\invar{H}(\Rder[q]\invar{N}(A)) \quad \Rightarrow \quad \Rder\invar{G}(A).
\end{equation*}
Utilizing this powerful tool, we obtain a proof of \Cref{thm:isom_gc} by replacing $A$ with $A^H$. In this case, the $E_2$ page is given by
\begin{equation*}
    \Rder\invar{H}(\Rder\invar{N}(A^H)) \cong \Rder\invar{H}(\Rder\invar{N}(A)^H) \cong \Rder\invar{N}(A).
\end{equation*}
The first equality is true since the $N$-action on $H$ is trivial, and the right derived functor $\Rder\invar{N}$ commutes with infinite product in any category of modules. The second equality is a consequence of \Cref{lem:gc_con_at_zero}. Since $E_2^{\bullet, \bullet}$ is concentrated at degrees $(0, \bullet)$, it follows that the Grothendieck spectral sequence degenerates at $E_2$ page, and that $ \Rder\invar{G}(A^H) \cong \Rder\invar{N}(A) $ canonically.

\subsection{Comparing the Group Cohomology of \texorpdfstring{$\Rper$}{R\_per} with the CKS Complex}

\begin{theorem}[thm:GC_CKS]{}
    For each $ n, k \in \mathbb{N} $, there is a non-canonical isomorphism
    \begin{equation*}
        \HH_{\mathrm{Grp}}^n(\HH_1(\Gamma), \Rper[2k]) \cong \bigoplus_{p+q=n+k} \HH^p(\grCKS[k][q], \dCKS).
    \end{equation*}
    \noskipline
    \begin{proof}
        By definition, $ \HH_{\mathrm{Grp}}^n(\HH_1(\Gamma), \Rper[2k]) = \Rder[n]\invar{\HH_1(\Gamma)}(\Rper[2k]) $. By \Cref{thm:HT_eqv}, there is a homotopy equivalence $ (\grHT[k][\bullet+k], \dHT) \to \Rper[2k] $, which is $\HH_1(\Gamma)$-equivariant by construction. By functoriality, this induces an isomorphism $ \totRder\invar{\HH_1(\Gamma)}(\grHT[k][\bullet+k]) \to \totRder\invar{\HH_1(\Gamma)}(\Rper[2k]) $ in the derived category. This means that given any projective resolution $P_\bullet$ of $\ZZ$ as a $\ZZ[G]$-module, the chain map $ \Hom_{\ZZ[G]}^\bullet(P_\bullet, \grHT[k][\bullet+k]) \to \Hom_{\ZZ[G]}^\bullet(P_\bullet, \Rper[2k]) $, where the domain is given the total differential, is a quasi-isomorphism.
        Consequently, for each $ k \in \mathbb{N} $, there is a spectral sequence of the double complex $\totRder\invar{\HH_1(\Gamma)}(\grHT[k][\bullet+k])$ given by
        \begin{equation*}
            E_2^{p-k, q} = \HH^{p}(\Rder[q]\invar{\HH_1(\Gamma)}(\grHT[k]), \dHT) \quad \Rightarrow \quad \Rder[\bullet]\invar{\HH_1(\Gamma)}(\Rper[2k]).
        \end{equation*}

        We claim that there is an isomorphism of complexes $ (\Rder[q]\invar{\HH_1(\Gamma)}(\grHT[k]), \dHT) \cong (\grCKS[k][q], \dCKS) $. Indeed,
        \begin{equation*}
            \Rder[q]\invar{\HH_1(\Gamma)}(\grHT[k]) = \Rder[q]\invar{\HH_1(\Gamma)}\left(\bigoplus_{S \in \face} \ZZ^{\ZZ^S} \otimes \Ext[k-\bullet] \HH_1(\Gamma \del S)\right) \cong \bigoplus_{S \in \face} \Rder[q]\invar{\HH_1(\Gamma)}(\ZZ^{\ZZ^S}) \otimes \Ext[k-\bullet] \HH_1(\Gamma \del S).
        \end{equation*}
        But the $\HH_1(\Gamma)$-action on $\ZZ^{\ZZ^S}$ is induced by the isomorphism $ \ZZ^S \cong \HH_1(\Gamma) / \HH_1(\Gamma \del S) $. Thus, \Cref{thm:isom_gc} implies that $ \Rder\invar{\HH_1(\Gamma)}(\ZZ^{\ZZ^S}) \cong \Rder\invar{\HH_1(\Gamma \del S)}(\ZZ) $, which evaluates to $ \Ext \HH^1(\Gamma \del S) $. The differentials coincide by construction. Thus, $ E_2^{p-k, q} \cong \HH^{p}(\grCKS[k][q], \dCKS) $ canonically.

        Note that $\grHT[k]$ is a product of free abelian groups by definition, hence flat. By \Cref{thm:degen_spectral_seq}, this spectral sequence degenerates at $E_2$ page and we have the desired isomorphism.
    \end{proof}
\end{theorem}

\begin{corollary}[cor:]{}
    The group cohomology $\HH_{\mathrm{Grp}}^\bullet(\HH_1(\Gamma), \Rper)$ is isomorphic to the integral cohomology ring of $\Mmul$.
\end{corollary}

We now compute an example of $\HH_{\mathrm{Grp}}^\bullet(\HH_1(\Gamma), \Rper)$.

\begin{example}[exp:gc_loopgraph]{}
    For $ \Gamma = \loopgraph $, we have
    \begin{equation*}
        \Rper = \frac{\widehat{\ZZ[\{e_i \mid i \in \ZZ\}]}}{\{ e_i e_j \mid i, j \in \ZZ, i \neq j \}, \sum_{i \in \ZZ} e_i}
    \end{equation*}
    and the action of $\HH_1(\Gamma) \cong \ZZ$ is given by $n \cdot e_i = e_{n+i}$. Note that we have $ e_i^2 = - e_i \sum_{j \in \ZZ, j \neq i} e_j = 0 $. Thus, $ \Rper[0] \cong \ZZ $, $ \Rper[>2] \cong 0 $ and $\Rper[2]$ fits into the $\ZZ$-equivariant short exact sequence
    \begin{equation*}
        \begin{tikzcd}
            0 \arrow[r] & \ZZ \arrow[r] & \ZZ^\ZZ \arrow[r] & \Rper[2] \arrow[r] & 0
        \end{tikzcd}
    \end{equation*}
    where the first $\ZZ$ has trivial action and includes as the diagonal, and the projection is given by $(a_i)_{i \in \ZZ} \mapsto \sum_{i \in \ZZ} a_i e_i$. Note that $\HH^\bullet_{\mathrm{Grp}}(\ZZ, \ZZ) \cong \Ext \ZZ$, and we have $\HH^\bullet_{\mathrm{Grp}}(\ZZ, \ZZ^\ZZ) \cong \ZZ$ by \Cref{lem:gc_con_at_zero}. It follows by the induced long exact sequence in group cohomology that $ \HH^\bullet_{\mathrm{Grp}}(\ZZ, \Rper[2]) \cong \ZZ $ is concentrated at degree $0$, which is in fact generated by the invariant element $ \sum_{i \in \ZZ} i e_i \in \Rper[2] $. It follows that $ \HH^0_{\mathrm{Grp}}(\ZZ, \Rper[0]) \cong \HH^1_{\mathrm{Grp}}(\ZZ, \Rper[0]) \cong \HH^0_{\mathrm{Grp}}(\ZZ, \Rper[2]) \cong \ZZ $, which coincides with the cohomology of the CKS complex as computed in \Cref{exp:cks_loopgraph}.
\end{example}

\subsection{Deletion-Contraction for Multiplicative Hypertoric Varieties}

Let us examine the action of $\HH_1(\Gamma)$ on $\Rper[\bullet][\Gamma\del e]^\ZZ$ in \Cref{thm:Rper_delcon}. Writing $ (r_i)_{i \in \ZZ} \in \Rper[\bullet][\Gamma\del e]^\ZZ $ by $ \sum r_i 1_{e_i} $, we have
\begin{equation*}
    \gamma \cdot \sum r_i 1_{e_i} = \sum \gamma \cdot r_i 1_{\gamma \cdot e_i} = \sum \gamma \cdot r_i 1_{e_{i - \pair{\gamma}{e}}}.
\end{equation*}
Hence, the $\HH_1(\Gamma)$-action on $\ZZ$ is by identifying $ \HH_1(\Gamma) / \HH_1(\Gamma \del e) \cong \ZZ $, as in \Cref{eqn:del_seq}. By \Cref{thm:isom_gc}, we have
\begin{equation*}
    \Rder\invar{\HH_1(\Gamma)}(\Rper[\bullet][\Gamma\del e]^\ZZ) \cong \Rder\invar{\HH_1(\Gamma \del e)}(\Rper[\bullet][\Gamma\del e]).
\end{equation*}
Thus, applying $\Rder\invar{\HH_1(\Gamma)}$ to \Cref{thm:Rper_delcon} gives the following long exact sequence.

\begin{theorem}[thm:]{}
    If $e\in E$ is neither a bond nor a loop, we have a deletion-contraction long exact sequence.
    \begin{equation*}
        \begin{tikzcd}
            \cdots \arrow[r] & \Rder\invar{\HH_1(\Gamma \del e)}(\Rper[\bullet-2][\Gamma\del e]) \arrow[r] & \Rder\invar{\HH_1(\Gamma)}(\Rper) \arrow[r] & \Rder\invar{\HH_1(\Gamma \con e)}(\Rper[\bullet][\Gamma\con e]) \arrow[r] & \cdots
        \end{tikzcd}
    \end{equation*}
\end{theorem}

We suspect that this deletion-contraction sequence is related to the deletion-contraction sequence for $\CKS$ in \Cref{thm:cks_del_con} under some isomorphisms as in \Cref{thm:GC_CKS}.
\minormissing{compare GC and CKS del-con}